\documentclass[12pt,oneside,reqno]{amsart}
\usepackage{}
\usepackage{amssymb}
\usepackage{bbm}
\usepackage{cases}
\usepackage{amsmath}
\usepackage{graphicx}
\usepackage{mathrsfs}
\usepackage{stmaryrd}
\usepackage{amsfonts}
\usepackage{enumerate,amsmath,amssymb,amsthm}
\usepackage[numbers,sort&compress]{natbib}

\numberwithin{equation}{section}

\newcommand{\be}{\begin{eqnarray}}
\newcommand{\ee}{\end{eqnarray}}
\newcommand{\ce}{\begin{eqnarray*}}
\newcommand{\de}{\end{eqnarray*}}
\newtheorem{theorem}{Theorem}[section]
\newtheorem{lemma}[theorem]{Lemma}
\newtheorem{remark}[theorem]{Remark}
\newtheorem{definition}[theorem]{Definition}
\newtheorem{proposition}[theorem]{Proposition}
\newtheorem{example}[theorem]{Example}
\newtheorem{corollary}[theorem]{Corollary}

\def\e{{\mathrm{e}}}
\def\eps{\varepsilon}

\def\p{\partial}

\def\[{{\Big[}}
\def\]{{\Big]}}
\def\<{{\langle}}
\def\>{{\rangle}}
\def\({{\Big(}}
\def\){{\Big)}}

\def\bx{{\mathbf{x}}}

\def\dif{{\mathord{{\rm d}}}}

\def\no{\nonumber}
\def\={&\!\!=\!\!&}
\def\bt{\begin{theorem}}
\def\et{\end{theorem}}
\def\bl{\begin{lemma}}
\def\el{\end{lemma}}
\def\br{\begin{remark}}
\def\er{\end{remark}}

\def\bd{\begin{definition}}
\def\ed{\end{definition}}
\def\bp{\begin{proposition}}
\def\ep{\end{proposition}}
\def\bc{\begin{corollary}}
\def\ec{\end{corollary}}
\def\bx{\begin{example}}
\def\ex{\end{example}}

\def\cB{{\mathcal B}}

\def\cI{{\mathcal I}}
\def\cJ{{\mathcal J}}

\def\cM{{\mathcal M}}
\def\cN{{\mathcal N}}

\def\cT{{\mathcal T}}

\def\mE{{\mathbb E}}

\def\mH{{\mathbb H}}
\def\mI{{\mathbb I}}

\def\mN{{\mathbb N}}

\def\mP{{\mathbb P}}

\def\mR{{\mathbb R}}

\def\mW{{\mathbb W}}

\def\sC{{\mathscr C}}
\def\sD{{\mathscr D}}
\def\sE{{\mathscr E}}
\def\sF{{\mathscr F}}

\def\sK{{\mathscr K}}
\def\sL{{\mathscr L}}

\def\sS{{\mathscr S}}
\def\sT{{\mathscr T}}

\def\geq{\geqslant}
\def\leq{\leqslant}

\allowdisplaybreaks

\begin{document}

\title{Irregular Stochastic differential equations driven by a family of Markov processes}

\date{}

\author{{Longjie Xie}
and  {Lihu Xu}
}

\address{Longjie Xie:
School of Mathematics and Statistics, Jiangsu Normal University,
Xuzhou, Jiangsu 221000, P.R.China\\
Email: xlj.98@whu.edu.cn
 }
\address{Lihu Xu: 1. Department of Mathematics, Faculty of Science and Technology,
University of Macau, Av. Padre Tom\'{a}s Pereira, Taipa Macau, China. 2. UM Zhuhai Research Institute, Zhuhai, 519080, China \\
Email: lihuxu@umac.mo, xulihu2007@gmail.com}

\thanks{The first author is supported by the Project
Funded by PAPD of Jiangsu Higher Education Institutions. The second author is partially supported by the grants NSFC No 11571390, Macao S.A.R. FDCT/030/2016/A1 and University of Macau MYRG2016-00025-FST}

\begin{abstract}
Using heat kernel estimates, we prove the pathwise uniqueness for strong solutions of irregular stochastic differential equation driven by a family of Markov process, whose generator is a non-local and non-symmetric L\'evy type operator. Due to the extra term $1_{[0,\sigma(X_{s-},z)]}(r)$ in multiplicative noise, we need to derive some new regularity results for the generator and use a trick of mixing $L_1$ and $L_2$-estimates by Kurtz and Protter \cite{Ku-Po}.

  \bigskip


\end{abstract}

\maketitle \rm

\section{Introduction}

Nowadays, much attentions have been paid to the non-local operators and their corresponding pure jump processes, as these processes are more realistic models for many practice applications. Consider the following non-local and non-symmetric L\'evy type operator:
\begin{align}
\sL f(x):=\sL_\nu^\sigma f(x)+b(x)\cdot\nabla f(x),\quad\forall f\in C^\infty_0(\mR^d),  \label{oper}
\end{align}
where $b(x)$ is a measurable function and
\begin{align*}
\sL^\sigma_{\nu} f(x):=\int_{\mR^d}\big[f(x+z)-f(x)-1_{\{|z|\leq 1\}}z\cdot\nabla f(x)\big]\sigma(x,z)\nu(\dif z).
\end{align*}
Here, $\nu$ is a L\'evy measure on $\mR^d$ satisfying
$$
\int_{\mR^d\setminus \{0\}}(|z|^2\wedge 1)\nu(\dif z)<\infty,
$$
and $\sigma:\mR^d\times\mR^d\rightarrow\mR$ is measurable.
The operator $\sL$ is a non-local version of the classical second order elliptic operator with non-divergence form and has been intensely studied in the last decade by people in the community of analysis and PDEs, see \cite{E-I-K,Ko} and the references therein.
While from the probability point of view, it is known via the martingale method (see\cite{M-P}) that under certain assumptions on $\nu, \sigma$ and $b$, there exists a Markov process $X_t$ with $\sL$ as its generator, and the measure $\nu$ describes the jumps of $X_t$. It is natural to ask  wether one can construct $X_t$ via the It\^o's calculus so that we can have another look at $\sL$ from the view of stochastic differential equations (SDEs). However, the classical SDE driven by pure jump L\'evy process is not very suitable (see more discussions in \cite[Section 1]{Xie}). Its connection to SDE was found very recently.

To specify the SDE that we are going to study, let
$\cN(\dif z,\dif r,\dif t)$ be a Poisson random measure on $\mR^d\times[0,\infty)\times[0,\infty)$ with intensity measure $\nu(\dif z)\dif r\dif t$, and
$\tilde \cN(\dif z,\dif r,\dif t):=\cN(\dif z,\dif r,\dif t)-\nu(\dif z)\dif r\dif t$ is the compensated Poisson random measure.
Then, the Markov process $X_t$ corresponding to $\sL$ should satisfy the following SDE:
\begin{align}
\dif X_t&=\int_0^{\infty}\!\!\!\!\int_{|z|\leq 1}\!\!1_{[0,\sigma(X_{s-},z)]}(r)z\tilde \cN(\dif z,\dif r,\dif t)\no\\
&\quad+\int_0^{\infty}\!\!\!\!\int_{|z|> 1}\!\!1_{[0,\sigma(X_{s-},z)]}(r)z \cN(\dif z,\dif r,\dif t)+b(X_t)\dif t, \ \ \ \ X_0=x\in\mR^d. \label{sde2}
\end{align}
In fact, noticing that for a function $f$ and any $r>0$,
\begin{align*}
f\big(x+1_{[0,\sigma(x,z)]}(r)z\big)-f(x)=1_{[0,\sigma(x,z)]}(r)\big[f(x+z)-f(x)\big],
\end{align*}
an application of It\^o's formula shows that the generator of the solution to SDE (\ref{sde2}) is given exactly by (\ref{oper}). Note that the driven noise is a Markov process but not necessarily L\'evy type \cite{Kurz2}. This makes such kind of SDEs more interesting and are worthy of study.
Under the conditions that $b$ is bounded and global Lipschitz continuous, and $\sigma$ is bounded with
\begin{align}
\int_{\mR^d}|\sigma(x,z)-\sigma(y,z)|\cdot|z|\nu(\dif z)\leq C_1|x-y|, \ \ \ \ \forall \ \ x, y \in \mR^d,  \label{kur}
\end{align}
and some other assumptions, Kurtz \cite{Kurz2} showed the existence and uniqueness of strong solution to SDE \eqref{sde2}, see also \cite{Im-Wi,Ku-Po} for related results and applications.

\vskip 3mm

Our aim in this paper is to prove that SDE \eqref{sde2} admits a unique strong solution under some weak assumptions on the coefficients $\sigma$ and $b$ as well as the jump measure $\nu$, from which we can see the regularization effect of such kind of noises on the deterministic systems.

The irregular SDEs driven by pure jump noises have been extensively studied in the past several decades.
Note that when $d=1$ and $L_t$ is a symmetric $\alpha$-stable process with $\alpha<1$, Tanaka, Tsuchiya and Watanabe \cite{Ta-Ts-Wa} showed that if $b$ is bounded and $\beta$-H\"older continuous with $\alpha+\beta<1$, SDE
\begin{align}
\dif X_t=\dif L_t+b(X_t)\dif t,\quad X_0=x\in\mR^d \label{levy}
\end{align}
may not have pathwise uniqueness strong solutions. When $\alpha\geq 1$, $b$ is bounded and $\beta$-H\"older continuous with
$\beta>1-\frac{\alpha}{2},$
it was proved by Priola \cite{Pri} that there exists a unique strong solution $X_t(x)$ to SDE (\ref{levy}) for each $x\in\mR^d$. Recently, Zhang \cite{Zh00} obtained the pathwise uniqueness to SDE (\ref{levy}) when $\alpha>1$ and $b$ is a bounded function in some local Sobolev space. See also \cite{B-B-C,Ch-So-Zh,Pri2} for related results. We also would like to mention the paper \cite{M-X} where SDEs driven by multiplicative L\'evy noise with Lipschitz diffusion coefficient and H\"older drift was considered.
For the study of irregular SDEs driven by Brownian motion, we refer readers to \cite{Fa-Lu-Th,Fe-Fl-2,Fl-Gu-Pr,Kr-Ro,M-N-P-Z, W,Wa,XZ,Zh3,Zh1}.

Let us compare our results with the literatures above. To prove the uniqueness of strong solution, we shall follow a well known strategy \cite{Ch-So-Zh,Kr-Ro, Pri, Zh1}, which is to derive a new SDE with better coefficients by Zvonkin's transformation and get the uniqueness of the original equation from the new one. The crucial point of this approach is to study the regularity of transformation equations which vary with different SDEs. There are several new aspects that we would like to stress for SDE \eqref{sde2} as the following.

First of all, our main tool for studying the transformation equation (see \eqref{pide1} below) is the heat kernel (also called fundamental solution) of the operator $\sL_\nu^\sigma$, it seems the first time to use heat kernel estimates to study the pathwise uniqueness irregular SDEs, see \cite{C-K-K,K-S} for the study of weak uniqueness of SDEs with L\'evy noise by using heat kernels. Secondly,  all the above works are for singular SDEs driven by Brownian motions or additive L\'evy noises, in the latter case, one only needs to study the operator $\sL_0$ defined by
$$
\sL_0 f(x):=\int_{\mR^d}\Big[f(x+z)-f(x)-1_{\{|z|\leq 1\}}z\cdot\nabla f(x)\Big]\nu(\dif z),\quad\forall f\in C_0^\infty(\mR^d).
$$
The analysis relies on
the nice property of $\sL_0$ and the $C^2$ smoothing effect of its semigroup. However, we study the multiplicative noise and the semigroup generated by $\sL_{\nu}^\sigma$ only has $C^{\alpha+\beta}$ regularity with $\alpha+\beta<2$ (see Remark \ref{ree}), we need to use more delicate analysis and interpolation theorems to fit our less regularity property into the frame of  Zvonkin's argument. We mention that in \cite{Xie}, the first author consider the same SDE with critical case $\alpha=1$ and $b$ in H\"older spaces, here we shall
consider $\alpha\in(1,2)$ but with a more irregular drift term $b$ in fractional Sobolev spaces, and the proof in this paper is more involved.
Lastly, when proving the Krylov-type estimate and performing Zvonkin transformation, we needs to solve
a semi-linear elliptic equation and the resolvent equation of $\sL_\nu^\sigma$
in the framework of Sobolev
space. Because a well developed elliptic equation theory as in \cite{Kr-Ro,Zh3,Zh1} is not available
for $\sL_\nu^\sigma$, we derive a generalized It\^o's formula for H\"older functions and solve the
corresponding integral equation in Sobolev spaces.

Another novelty in our analysis is the technique for handling the extra term $1_{[0,\sigma(X_{s-},z)]}(r)$ when we prove the pathwise uniqueness in the last section. The usual $L_2$-estimate in the known literatures is not applicable. Fortunately we can use a trick of mixing $L_1$ and $L_2$-estimate as the replacement \cite{Kurz2,Ku-Po}. Due to the irregularity of $b$ and $\sigma$, it is much more complicated than \cite{Kurz2} to apply this trick.

 \vspace{2mm}
Finally, we mention that studying the unique strong solution of SDE \eqref{sde2} with irregular coefficients not only has its own interests but also helps to better understand the nonlocal operator $\sL$ (\cite{Im-Wi}). Another motivation for studying SDE \eqref{sde2} is because of the special noise. As mentioned above, the driven noise is a Markov process but not necessarily L\'evy type. This has been found very useful in applications, for instance, Markov type noise plays a crucial role as the control when proving Freidlin-Wentzell type large deviation for L\'evy type SDEs via weak convergence approach \cite{BDM11, BCD13, ZhZh15}.

 \vspace{2mm}

The organization of the paper is as the following. Section 2 gives the main result with some comments and comparisons with known literatures. Sections 3 and 4 are both preparation sections, the former for some estimates of heat kernel of $\sL_\nu^\sigma$ and the latter for the regularity of the corresponding semigroup $\cT_t$.  Krylov's estimate and Zvonkin's transformation are studied in the 5th section and applied in the last section to prove the strong uniqueness of SDE \eqref{sde2}.   Throughout this paper, we use the following convention: $C$ with or without subscripts will denote a positive constant, whose value may change in different places, and whose dependence on parameters can be traced from calculations.

\vskip 2mm
{\bf Acknowledgements:} We would like to gratefully thank Professors Rengming Song and Xicheng Zhang for very helpful discussions.

\section{Main result}
We assume that for all $x\in\mR^d$,
\begin{align}
\sigma(x,z)=\sigma(x,-z),\quad\forall z\in\mR^d, \label{s2}
\end{align}
and that there exists a function $\tilde\kappa$ such that
\begin{align}
\nu(\dif z)=\frac{\tilde\kappa(z)}{|z|^{d+\alpha}}\dif z, \quad \tilde\kappa(z)=\tilde\kappa(-z),\quad \kappa_0\leq \tilde\kappa(z)\leq \kappa_1,\label{nu}
\end{align}
with $\alpha\in(1,2)$ and $\kappa_0, \kappa_1$ are two positive constants. The symmetric in $z$ of $\sigma$ and $\tilde\kappa$ is a common assumption in the literature, see \cite{Ca-Si}. As a result, we can also write $\sL_\nu^\sigma$ as
\begin{align}
\sL^\kappa_{\alpha}\varphi(x)&=\text{p.v.}\int_{\mR^d}\big[\varphi(x+z)-\varphi(x)\big]\sigma(x,z)\nu(\dif z)\no\\
&=\text{p.v.}\int_{\mR^d}\big[\varphi(x+z)-\varphi(x)\big]\frac{\kappa(x,z)}{|z|^{d+\alpha}}\dif z,\label{op}
\end{align}
where
\begin{align}
\kappa(x,z)=\sigma(x,z)\tilde\kappa(z). \label{kappa}
\end{align}
The operator $\sL^\kappa_{\alpha}$ is a non-local and non-symmetric operator, which can be seen as a generalization of variable coefficients fractional Laplacian operator.

\vspace{3mm}
For brevity, we set $B_n:=\{x\in\mR^d:|x|\leq n\}$. Our main result is:

\bt\label{main}
Let (\ref{s2}) hold and the L\'evy measure $\nu$ satisfies (\ref{nu}). Suppose that for any $n\in\mN$:

\begin{enumerate}
\item [{\bf (H$\sigma$)}] There exists a function $\zeta\in L^{q}(B_n)$ with $q>d/\alpha$, such that for almost all $x,y\in B_n$,
\begin{align}
\int_{\mR^d}|\sigma(x,z)-\sigma(y,z)|(|z|\wedge1)\nu(\dif z)\leq |x-y|\Big(\zeta(x)+\zeta(y)\Big),\label{a1}
\end{align}
and for some constants $k^n_0, k^n_1>0$, $\beta\in(0,1)$ and $C_n>0$,
\begin{align}
k^n_0\leq \sigma(x,z)\leq k^n_1,\quad|\sigma(x,z)-\sigma(y,z)|\leq C_n|x-y|^{\beta},\quad \forall x,y\in B_n,\,\,\forall z\in\mR^d.\label{ho}
\end{align}

\item [{\bf (Hb)}] For some $\theta\in (1-\frac{\alpha}{2},1)$ and $p>2d/\alpha$,
\begin{align*}
\int_{B_n}\!\int_{B_n}\frac{|b(x)-b(y)|^p}{|x-y|^{d+\theta p}}\dif x\dif y<+\infty,
\end{align*}
and it holds
\begin{align*}
\sup_{x\in B_n}|b(x)|<\infty.
\end{align*}
\end{enumerate}
\vspace{-2mm}
Then, for each $x\in\mR^d$, there exists an stopping time $\varsigma(x)$ (called the explosion time) and a unique strong solution $X_t(x)$ to SDE (\ref{sde2}) such that
\begin{align}
\lim_{t\uparrow\varsigma(x)}X_t(x)=\infty,\quad a.s..    \label{xt}
\end{align}
\et
\vskip 3mm

Let us make some comments on the assumptions
and give an example for our result with a comparison with known literatures.
\br
(1). It is clear that the assumption (\ref{a1}) is a generalization of condition (\ref{kur}) in \cite{Kurz2}.
For an very interesting example of $\sigma$, we can take
$$
\sigma(x,z)=K(z)+\tilde\sigma(x)|z|^{\gamma} \ \ {\rm for} \ \  |z| \le 1, \ \ \ \ \ \sigma(x,z)=K(z)+\tilde\sigma(x) \ \ {\rm for} \ \  |z|>1,
$$
with $0<K_1\leq K(z)\leq K_2$, $\gamma>\alpha-1$ and $\nabla\tilde\sigma\in L^q_{loc}(\mR^d)$ with $q>d/\alpha$, where $\nabla$ denotes the weak derivative. Since we assume $\alpha>1$, our theorem can cover the regime $q \in (d/\alpha, d]$. However, for the following SDE driven by multiplicative Brownian motion \cite{Zh3}:
\begin{equation*}
\dif X_t=\sigma(X_t)\dif W_t+b(X_t)\dif t,\quad X_0=x\in\mR^d,
\end{equation*}
one has to assume that $\nabla\sigma\in L^q(\mR^d)$ with $q>d$. Here, the main trick is that $\sigma$ appears in the term $1_{[0,\sigma(X_{s-},z)]}(r)$. When $\nabla\tilde\sigma\in L^q_{loc}(\mR^d)$ with $q>d$, we can also have the H\"older continuity in (\ref{ho}) by the Sobolev embedding theorem.

(2). For interesting examples of irregular drift coefficient $b$, we can take $b(x)=1_A(x)$ for a measurable set $A \in \mR^d$, see \cite[Remark 1.2]{Zh00} for details.

(3).  The conditions (\ref{s2}), (\ref{nu}) and (\ref{ho}) are assumed so that we can use the results obtained in \cite{Ch-Zh}. Under (\ref{s2}), (\ref{nu}) and the global assumptions
\begin{align}
0<\tilde k_0\leq \sigma(x,z)\leq \tilde k_1,\,\,\,|\sigma(x,z)-\sigma(y,z)|\leq C_0|x-y|^{\beta},\quad \forall x,y\in \mR^d,\,\,\forall z\in\mR^d, \label{s1}
\end{align}
it was proved that there exists a unique fundamental solution $p(t,x,y)$ for $\sL^\kappa_\alpha$, see \cite[Theorem 1.1]{Ch-Zh}. Here, we only need the local boundness and the local H\"older continuity of $\sigma$ in (\ref{ho}) thanks to the stopping time technique. Furthermore, we shall prove better regularities of $p(t,x,y)$ (see Theorem \ref{seen}) than obtained in \cite{Ch-Zh},  which seem to be new and have independent interests.
\er

\section{Heat kernel estimates}

We briefly recall the construction of the heat kernel $p(t,x,y)$ for operator $\sL^\kappa_{\alpha}$ in \cite{Ch-Zh}, from which we derive some important estimates of $p(t,x,y)$ (see Theorem \ref{seen} below) for further use in next sections.
From now on, we assume that (\ref{s2}), (\ref{nu}) and (\ref{s1}) always hold.

First of all, in view of (\ref{s2}) (\ref{nu}) and (\ref{op}), we can also write
\begin{align}
\sL^\kappa_{\alpha} \varphi(x)=\frac{1}{2}\int_{\mR^d}\delta_\varphi(x,z)\frac{\kappa(x,z)}{|z|^{d+\alpha}}\dif z,\quad\forall \varphi\in C^\infty_0(\mR^d), \label{opp}
\end{align}
where
$$
\delta_\varphi(x,z):=\varphi(x+z)+\varphi(x-z)-2\varphi(x).
$$
In order to reflect the dependence of $\kappa$ with respect to $x$, we shall also use $\sL_\alpha^{\kappa,x}$ instead of $\sL_\alpha^\kappa$.
To shorten the notation, we set for $\gamma, \beta\in \mR$,
$$
\varrho_\gamma^\beta(t,x):=t^{\frac{\gamma}{\alpha}}\big(|x|^{\beta}\wedge 1\big)\big(|x|+t^{1/\alpha}\big)^{-d-\alpha}.
$$
The following 3-P type inequalities shall be used below from time to time.

\bl
(i). If $\gamma_1+\beta_1>0$ and $\gamma_2+\beta_2>0$, then there exists a constant $C_1>0$ such that for all $t\geq 0$ and $x,y\in\mR^d$,
\begin{align}
&\int_0^t\!\!\!\int_{\mR^d}\varrho_{\gamma_1}^{\beta_1}(t-s,x-z)\varrho_{\gamma_2}^{\beta_2}(s,z-y)\dif z\dif s\no\\
&\leq C_1\Big(\varrho^{0}_{\gamma_1+\gamma_2+\beta_1+\beta_2}+\varrho_{\gamma_1+\gamma_2+\beta_2}^{\beta_1}+ \varrho^{\beta_2}_{\gamma_1+\gamma_2+\beta_1}\Big)(t,x-y). \label{3p}
\end{align}
(ii). For all $\beta_1,\beta_2\in[0,\alpha]$ and $\gamma_1,\gamma_2\in\mR$, there exists a constant $C_2>0$ such that for any $t\geq 0$ and $x\in\mR^d$,
\begin{align}
&\int_{\mR^d}\varrho_{\gamma_1}^{\beta_1}(t,x-z)\varrho_{\gamma_2}^{\beta_2}(t,z)\dif z\no\\
&\leq C_2\Big(\varrho^{0}_{\gamma_1+\gamma_2+\beta_1+\beta_2-\alpha}+\varrho_{\gamma_1+\gamma_2+\beta_2-\alpha}^{\beta_1}+ \varrho^{\beta_2}_{\gamma_1+\gamma_2+\beta_1-\alpha}\Big)(t,x). \label{3p1}
\end{align}
\el
\begin{proof}
The first inequality is given by \cite[Lemma 2.1 (iii)]{Ch-Zh}, while the second one can be proved entirely by the same arguments as \cite[Lemma 2.1 (ii)]{Ch-Zh}, the details are omitted.
\end{proof}

Let $p_{\alpha}(t,x)$ denote the heat kernel of operator $\Delta^{\frac{\alpha}{2}}$ (or equivalently, the transition density of $d$-dimensional symmetric $\alpha$-stable process). It is well known that there exists a constant $C_0$ such that
\begin{align}
C_0^{-1}\varrho_\alpha^0(t,x)\leq p_{\alpha}(t,x)\leq C_0\varrho_\alpha^0(t,x), \label{k0}
\end{align}
and for every $k\in\mN$, it holds for some $C_k$ that
\begin{align}
|\nabla^k p_{\alpha}(t,x)|\leq C_k\varrho_{\alpha-k}^0(t,x). \label{kk}
\end{align}
Set for $z\in\mR^d$,
\begin{align*}
\delta_{p_\alpha}(t,x;z):=p_\alpha(t,x+z)+p_\alpha(t,x-z)-2p_\alpha(t,x).
\end{align*}
It was shown by \cite[Lemma 2.2]{Ch-Zh} that that there exists a constant $C>0$ such that
\begin{align}
|\delta_{p_\alpha}(t,x;z)|\leq C\Big((t^{-\frac{2}{\alpha}}|z|^2)\wedge 1\Big)\Big(\varrho^0_\alpha(t,x\pm z)+\varrho^0_\alpha(t,x)\Big).\label{de}
\end{align}
With this estimate in hand and following the same ideas as in the proof of \cite[Theorem 2.4]{Ch-Zh}, we can derive the fractional derivative estimate of $p_{\alpha}(t,x)$. For completeness, we sketch the details here.

\bl\label{fe}
For any $0<\gamma<2$, there exists a constant $C_\gamma$ such that
\begin{align}
|\Delta^{\frac{\gamma}{2}}p_{\alpha}(t,x)|\leq C_\gamma\varrho^0_{\alpha-\gamma}(t,x). \label{e1}
\end{align}
\el
\begin{proof}
We may assume that $t\leq 1$, since the general case follows by the Chapman-Kolmogorov equation. By the definition of fractional Laplacian and as (\ref{opp}), we can write for $0<\gamma<2$,
\begin{align}
\Delta^{\frac{\gamma}{2}}p_{\alpha}(t,x)=\frac{c_{d,\gamma}}{2}\int_{\mR^d}\delta_{p_\alpha}(t,x;z)\frac{1}{|z|^{d+\gamma}}\dif z.\label{22}
\end{align}
Consequently, we have by (\ref{de})
\begin{align*}
|\Delta^{\frac{\gamma}{2}}p_\alpha(t,x)|&\leq C_{d,\gamma}\int_{\mR^d}|\delta_{p_\alpha}(t,x;z)|\cdot|z|^{-d-\gamma}\dif z\\
&\leq C_1\varrho^0_\alpha(t,x)\int_{\mR^d}\Big((t^{-\frac{2}{\alpha}}|z|^2)\wedge 1\Big)|z|^{-d-\gamma}\dif z\\
&\quad+C_1\!\!\int_{\mR^d}\Big((t^{-\frac{2}{\alpha}}|z|^2)\wedge 1\Big)\varrho^0_\alpha(t,x\pm z)|z|^{-d-\gamma}\dif z=:I_1+I_2.
\end{align*}
For $I_1$, by the assumption that $\gamma<2$, one can check easily that
\begin{align*}
I_1&=C_1\varrho^0_{\alpha-2}(t,x)\!\int_{|z|\leq t^{1/\alpha}}|z|^{2-d-\gamma}\dif z+C_1\varrho^0_{\alpha}(t,x)\!\int_{|z|> t^{1/\alpha}}|z|^{-d-\gamma}\dif z\leq C_2\varrho^0_{\alpha-\gamma}(t,x).
\end{align*}
As for the second term, similarly we write
\begin{align*}
I_2&=C_1t^{-\frac{2}{\alpha}}\!\!\int_{|z|\leq t^{1/\alpha}}\varrho^0_\alpha(t,x\pm z)|z|^{2-d-\gamma}\dif z+C_1\!\!\int_{|z|> t^{1/\alpha}}\varrho^0_\alpha(t,x\pm z)|z|^{-d-\gamma}\dif z=:I_{21}+I_{22}.
\end{align*}
We further control $I_{21}$ by
\begin{align*}
I_{21}&\leq C_3\varrho^0_{\alpha-2}(t,x)\!\!\int_{|z|\leq t^{1/\alpha}}|z|^{2-d-\gamma}\dif z\leq C_4\varrho^0_{\alpha-\gamma}(t,x).
\end{align*}
For $I_{22}$, if $|x|\leq 2t^{1/\alpha}$, then
\begin{align*}
I_{22}&\leq C_3t^{-\frac{d}{\alpha}}\!\int_{|z|> t^{1/\alpha}}|z|^{-d-\gamma}\dif z\leq C_3t^{-\frac{d+\gamma}{\alpha}}\leq C_4\varrho^0_{\alpha-\gamma}(t,x).
\end{align*}
If $|x|> 2t^{1/\alpha}$, we can deduce that
\begin{align*}
I_{22}&\leq C_1\left(\int_{\frac{|x|}{2}>|z|> t^{1/\alpha}}+\int_{|z|> \frac{|x|}{2}}\right)\varrho^0_\alpha(t,x\pm z)|z|^{-d-\gamma}\dif z\\
&\leq C_2t\int_{\frac{|x|}{2}>|z|> t^{1/\alpha}}\big(|x\pm z|+t^{1/\alpha}\big)^{-d-\alpha}|z|^{-d-\gamma}\dif z+C_2|x|^{-d-\gamma}\int_{|z|> \frac{|x|}{2}} \varrho^0_\alpha(t,x\pm z)\dif z\\
&\leq C_3\varrho^0_{\alpha}(t,x)\int_{|z|> t^{1/\alpha}}|z|^{-d-\gamma}\dif z+C_3|x|^{-d-\gamma}\leq C_4\varrho^0_{\alpha-\gamma}(t,x).
\end{align*}
Combing the above computations, we get (\ref{e1}).
\end{proof}

Now we fix $y\in\mR^d$, consider the freezing operator
$$
\sL_{\alpha}^{\kappa,y}f(x):=\text{p.v.}\int_{\mR^d}[f(x+z)-f(x)]\frac{\kappa(y,z)}{|z|^{d+\alpha}}\dif z,
$$
where $\kappa$ is given by (\ref{kappa}). It is known that there exists a symmetric $\alpha$-stable like process corresponding to $\sL_{\alpha}^{\kappa,y}$. Let $p_y(t,x)$ be the heat kernel of operator $\sL_{\alpha}^{\kappa,y}$. Since $\kappa$ is uniformly bounded, it follows from \cite[Theorem 1.1]{Ch-Ku} that for some constant $C_0$ independent of $y$,
\begin{align}
C_0^{-1}\varrho_\alpha^0(t,x)\leq p_{y}(t,x)\leq C_0\varrho_\alpha^0(t,x). \label{p0}
\end{align}
Moreover, if we set
$$
\hat\kappa(y,z):=\kappa(y,z)-\frac{\tilde k_0\kappa_0}{2},
$$
where $\tilde k_0, \kappa_0$ are the constants in (\ref{s1}) and (\ref{nu}), respectively, and let $\hat p_y(t,x)$ be the heat kernel of operator $\sL_{\alpha}^{\hat\kappa,y}$, by the construction of L\'evy process, we can write
\begin{align}
p_y(t,x)=\int_{\mR^d}p_{\alpha}(\tfrac{\tilde k_0\kappa_0}{2}t,x-z)\hat p_y(t,z)\dif z, \label{p00}
\end{align}
see also \cite[(2.23)]{Ch-Zh}. The advantage of (\ref{p00}) is that we can derive certain estimates for $p_y(t,x)$ by using properties of $p_\alpha(t,x)$. As an easy result, we have the following fractional derivative estimate of $p_y(t,x)$ and the H\"older continuity of $\nabla p_y(t,x)$. Here and below, both operators $\Delta^{\frac{\gamma}{2}}$ and $\nabla$ are acted with respect to the variable $x$.

\bl
For any $0<\gamma<2$, it holds
\begin{align}
|\Delta^{\frac{\gamma}{2}}p_{y}(t,x)|\leq C_\gamma\varrho^0_{\alpha-\gamma}(t,x). \label{e2}
\end{align}
and for any $\vartheta\in(0,1)$, $t>0$ and all $x,x',y\in\mR^d$,
\begin{align}
|\nabla p_y(t,x)-\nabla p_y(t,x')|\leq C_\vartheta|x-x'|^\vartheta\varrho_{\alpha-1-\vartheta}^0(t,\tilde x). \label{p0h}
\end{align}
where $C_\gamma, C_\vartheta$ are positive constants independent of $y$, and $\tilde x$ is the one of the two points $x$ and $x'$ which is nearer to zero point.
\el
\begin{proof}
It is enough to prove the estimates with $t\in(0,1)$. The first assertion can be verified by using the Fubini's theorem, (\ref{k0}), (\ref{e1}), (\ref{p0}), (\ref{p00}) (\ref{3p1}) and easy computations.
As for the second inequality, without lose of generality, we may assume that $|x|\leq |x'|$. In view of (\ref{kk}), we know that when $|x-x'|\geq t^{\frac{1}{\alpha}}/2$,
\begin{align*}
|\nabla p_\alpha(t,x)-\nabla p_\alpha(t,x')|&\leq C_1|x-x'|^\vartheta\Big(\varrho_{\alpha-1-\vartheta}^0(t,x)+\varrho_{\alpha-1-\vartheta}^0(t,x')\Big)\\
&\leq C_1|x-x'|^\vartheta\varrho_{\alpha-1-\vartheta}^0(t,x).
\end{align*}
While for $|x-x'|< t^{\frac{1}{\alpha}}/2$, we have by the mean value theorem that for some $\varepsilon\in[0,1]$,
\begin{align*}
|\nabla p_\alpha(t,x)-\nabla p_\alpha(t,x')|&\leq C_2|x-x'|\varrho_{\alpha-2}^0\big(t,x+\varepsilon(x'-x)\big)\\
&\leq C_2|x-x'|\varrho_{\alpha-2}^0(t,x)\leq C_2|x-x'|^\vartheta\varrho_{\alpha-1-\vartheta}^0(t,x).
\end{align*}
The desired estimate (\ref{p0h}) follows by (\ref{p00}), (\ref{p0}) and (\ref{3p1}).
\end{proof}

Let $\beta$ be the H\"older index in (\ref{s1}). Below, we always suppose that $\alpha+\beta<2$. This is assumed just to simplify the proof and is in fact not an restriction at all. Indeed, since we also assumed that $\sigma$ is bounded, (\ref{s1}) still holds true for any $\beta'<\beta$. Hence, it is enough to study the pathwise uniqueness of SDE (\ref{sde2}) when $\beta<2-\alpha$. We show the following estimate.

\bl\label{as}
Under (\ref{s1}), we have for $\gamma\in(0,2-\alpha)$ and all $x\in\mR^d$,
\begin{align}
\bigg|\int_{\mR^d}\Delta^{\frac{\alpha+\gamma}{2}}p_y(t,x-y)\dif y\bigg|\leq C_{d,\alpha,\gamma}t^{\frac{\beta-\gamma}{\alpha}-1},\label{00}
\end{align}
and for any $\vartheta\in(0,1)$ and $x, x'\in\mR^d$,
\begin{align}
\bigg|\int_{\mR^d}\Big[\nabla p_y(t,x-y)-\nabla p_y(t,x'-y)\Big]\dif y\bigg|\leq C_{d,\vartheta}|x-x'|^\vartheta t^{\frac{\beta-\vartheta-1}{\alpha}},\label{000}
\end{align}
where $C_{d,\alpha,\gamma}, C_{d,\vartheta}$ are positive constants.
\el
\begin{proof}
Since $\hat p_y(t,x)$ is a density function of Markov process, we have for any $\xi\in\mR^d$,
$$
\int_{\mR^d}\hat p_\xi(t,x-y)\dif y=1.
$$
Combing this with (\ref{22}), (\ref{p00}) and using the Fubini's theorem, it is easily checked that
\begin{align*}
\int_{\mR^d}\!\!\int_{\mR^d}\delta_{p_\xi}(t,x-y;z)\frac{c_{d,\alpha,\gamma}}{|z|^{d+\alpha+\gamma}}\dif z\dif y&=\int_{\mR^d}\Delta^{\frac{\alpha+\gamma}{2}}p_\xi(t,x-y)\dif y\\
&=\int_{\mR^d}\!\!\int_{\mR^d}\Delta^{\frac{\alpha+\gamma}{2}}p_\alpha(\tfrac{\tilde k_0\kappa_0}{2}t,x-y-z)\hat p_{\xi}(t,z)\dif z\dif y\\
&\!\!\!\stackrel{y+z=\tilde z}{=}\int_{\mR^d}\Delta^{\frac{\alpha+\gamma}{2}}p_\alpha(\tfrac{\tilde k_0\kappa_0}{2}t,x-z)\dif z=0,\quad \forall \xi\in\mR^d.
\end{align*}
As a result, we can write
\begin{align}
\sE_1:=\int_{\mR^d}\Delta^{\frac{\alpha+\gamma}{2}}p_y(t,x-y)\dif y&=\int_{\mR^d}\!\!\int_{\mR^d}\delta_{p_y}(t,x-y;z)\frac{c_{d,\alpha,\gamma}}{|z|^{d+\alpha+\gamma}}\dif z\dif y\no\\
&=\int_{\mR^d}\!\!\int_{\mR^d}\big[\delta_{p_y}(t,x-y;z)-\delta_{p_\xi}(t,x-y;z)\big]\Big|_{\xi=x}\frac{c_{d,\alpha,\gamma}}{|z|^{d+\alpha+\gamma}}\dif z\dif y.\label{ii}
\end{align}
By the proof of \cite[Theorem 2.5]{Ch-Zh}, we know that for any $0<\gamma'<\alpha$, there exists a $C_{\gamma'}$ such that
\begin{align*}
\big[\delta_{p_y}(t,x-y;z)-\delta_{p_\xi}(t,x-y;z)\big]\Big|_{\xi=x}&\!\leq C_{\gamma'}\big(|x-y|^\beta\wedge 1\big)\Big((t^{-\frac{2}{\alpha}}|z|^2)\wedge 1\Big)\\
&\times\Big(\big(\varrho^0_\alpha+\varrho^{\gamma'}_{\alpha-\gamma'}\big)(t,x-y\pm z)+\big(\varrho^0_\alpha+\varrho^{\gamma'}_{\alpha-\gamma'}\big)(t,x-y)\Big).
\end{align*}
Taking this into (\ref{ii}), choosing $\gamma'$ such that $\alpha+\gamma+\gamma'<2$ and arguing entirely the same as in the proof of Lemma \ref{fe}, we find that
\begin{align*}
\sE_1\leq C_{d,\alpha,\gamma}\!\int_{\mR^d}\varrho^{\beta}_{-\gamma}(t,x-y)\dif y\leq C_{d,\alpha,\gamma}t^{\frac{\beta-\gamma}{\alpha}-1}.
\end{align*}
Hence, (\ref{00}) is true. We proceed to prove (\ref{000}). Using (\ref{p00}) again, we write
\begin{align*}
\nabla p_y(t,x-y)-\nabla p_{y}(t,x'-y)=\int_{\mR^d}\sK_{\nabla p_{\alpha}}(t;x,x';y,z)\hat p_y(t,z)\dif z,
\end{align*}
where
$$
\sK_{\nabla p_{\alpha}}(t;x,x';y,z):=\nabla p_{\alpha}\big(\tfrac{\tilde k_0\kappa_0}{2}t,x-y-z\big)-\nabla p_{\alpha}\big(\tfrac{\tilde k_0\kappa_0}{2}t,x'-y-z\big).
$$
Let $\tilde x$ be the one of the two points $x$ and $x'$ which is nearer to $y+z$. Then, we know from the proof of (\ref{p0h}) that for any $\vartheta\in(0,1)$,
$$
|\sK_{\nabla p_{\alpha}}(t;x,x';y,z)|\leq C_\vartheta|x-x'|^\vartheta\varrho_{\alpha-1-\vartheta}^0\big(t,\tilde x-y-z\big).
$$
We may argue as in (\ref{ii}) to deduce that
\begin{align*}
\sE_2&:=\int_{\mR^d}\Big[\nabla p_y(t,x-y)-\nabla p_y(t,x'-y)\Big]\dif y\\
&=\int_{\mR^d}\!\int_{\mR^d}\sK_{\nabla p_{\alpha}}(t;x,x';y,z)\big[\hat p_y(t,z)-\hat p_{\xi}(t,z)\big]\Big|_{\xi=\tilde x}\dif z\dif y.
\end{align*}
Thanks to \cite[Theorem 2.5]{Ch-Zh}, we know that for any $0<\gamma'<\alpha$, there exists a $C_{\gamma'}$ such that
\begin{align*}
\big[\hat p_y(t,z)-\hat p_{\xi}(t,z)\big]\Big|_{\xi_=\tilde x}\leq C_{\gamma'}\big(|\tilde x-y|^{\beta}\wedge1\big)\Big(\varrho^{0}_{\alpha}+\varrho^{\gamma'}_{\alpha-\gamma'}\Big)(t,z),
\end{align*}
which yields by (\ref{3p1}) that
\begin{align*}
\sE_2&\leq C_\vartheta|x-x'|^\vartheta\int_{\mR^d}\!\int_{\mR^d}\varrho^0_{\alpha-1-\vartheta}(t,\tilde x-y-z)\Big(\varrho^{0}_{\alpha}+\varrho^{\gamma'}_{\alpha-\gamma'}\Big)(t,z)\dif z\cdot\big(|\tilde x-y|^{\beta}\wedge1\big)\dif y\\
&\leq C_\vartheta|x-x'|^\vartheta\int_{\mR^d}\Big(\varrho^{\beta}_{\alpha-1-\vartheta}+\varrho^{\gamma'+\beta}_{\alpha-1-\vartheta-\gamma'}\Big)\big(t,\tilde x-y\big)\dif y\leq C_\vartheta|x-x'|^\vartheta t^{\frac{\beta-\vartheta-1}{\alpha}}.
\end{align*}
The proof is finished.
\end{proof}

Now, the Levi's parametrix method suggests that the fundamental solution $p(t,x,y)$ of $\sL^{\kappa,x}_{\alpha}$ should be of the form
\begin{align}
p(t,x,y)=p_0(t,x,y)+\int_0^t\!\!\!\int_{\mR^d}p_0(t-s,x,z)q(s,z,y)\dif z\dif s, \label{heat}
\end{align}
where $p_0(t,x,y):=p_y(t,x-y)$ and $q(t,x,y)$ satisfies the integral equation
$$
q(t,x,y)=q_0(t,x,y)+\int_0^t\!\!\!\int_{\mR^d}q_0(t-s,x,z)q(s,z,y)\dif z\dif s
$$
with
$$
q_0(t,x,y):=\big(\sL^{\kappa,x}_{\alpha}-\sL_{\alpha}^{\kappa,y}\big)p_0(t,x,y).
$$
\vskip 3mm

The following lemma collects some estimates that we shall use below, whose proof can be found in \cite{Ch-Zh}.
\bl
The following statements hold:
\begin{enumerate}[(1).]
\item (\cite[Theorem 3.1]{Ch-Zh}) There exist constants $C_1, C_2$ such that for all $t\geq 0$ and $x,y\in\mR^d$,
\begin{align}
|q(t,x,y)|\leq C_1\big(\varrho_0^{\beta}+\varrho_\beta^0\big)(t,x-y),   \label{q1}
\end{align}
and for any $\gamma<\beta$, $t\geq 0$ and every $x,x',y\in\mR^d$,
\begin{align}
|q(t,x,y)-q(t,x',y)|\leq C_2\Big(|x-x'|^{\beta-\gamma}\wedge 1\Big)\Big(\big(\varrho^0_{\gamma}&+\varrho^{\beta}_{\gamma-\beta}\big)(t,x-y)\no\\
&\quad+\big(\varrho^0_{\gamma}+\varrho^{\beta}_{\gamma-\beta}\big)(t,x'-y)\Big),   \label{q2}
\end{align}
where $\beta$ is the constant in (\ref{ho}).

\item (\cite[Theorem 1.1]{Ch-Zh}) It hold for all $t>0$ and $x\in\mR^d$,
\begin{align}
\int_{\mR^d}p(t,x,y)\dif y=1, \label{11}
\end{align}
and there exists a constant $C_3$ such that
\begin{align}
|\nabla p(t,x,y)|\leq C_3\varrho^0_{\alpha-1}(t,x-y). \label{na}
\end{align}
\end{enumerate}
\el

In \cite{Ch-Zh}, it was also shown that for a constant $C>0$,
$$
|\Delta^{\frac{\alpha}{2}}p(t,x,y)|\leq C\varrho^0_{0}(t,x-y),
$$
in which the main point is to handle the singularity caused by the integral with respect to $s$. To study the strong solution of equation \eqref{sde2}, we need to prove more delicate estimates \eqref{es} and (\ref{ph}) as below, the proof is quite involved.

\bt\label{seen}
Suppose (\ref{s1}) holds true. Then, there exist constants $C_{d,\alpha,\gamma}, C_\vartheta>0$ such that for any $0\leq \gamma<\beta$,
\begin{align}
|\Delta^{\frac{\alpha+\gamma}{2}}p(t,x,y)|\leq C_{d,\alpha,\gamma}\varrho^0_{-\gamma}(t,x-y).   \label{es}
\end{align}
and for any $\vartheta\in(0,\alpha+\beta-1)$, $t>0$ and all $x,x',y\in\mR^d$,
\begin{align}
|\nabla p(t,x,y)-\nabla p(t,x',y)|\leq C_\vartheta|x-x'|^\vartheta\varrho_{\alpha-1-\vartheta}^0(t,\tilde x-y), \label{ph}
\end{align}
where $\tilde x$ is the one of the two points $x$ and $x'$ which is nearer to $y$.
\et
\begin{proof}
Still, we only consider the case when $t\leq 1$. For brevity, we set
$$
\sS(t,x,y):=\int_0^t\!\!\!\int_{\mR^d}p_0(t-s,x,z)q(s,z,y)\dif z\dif s.
$$
By Fubini's theorem, we can write
\begin{align*}
\big|\Delta^{\frac{\alpha+\gamma}{2}}\sS(t,x,y)\big|&\leq \Bigg|\!\int_{\frac{t}{2}}^t\!\!\!\int_{\mR^d}\Delta^{\frac{\alpha+\gamma}{2}}p_0(t-s,x,z)\Big(q(s,z,y)-q(s,x,y)\Big)\dif z\dif s\Bigg|\\
&\quad+\Bigg|\!\int_{\frac{t}{2}}^t\!\!\!\int_{\mR^d}\Delta^{\frac{\alpha+\gamma}{2}}p_0(t-s,x,z)\dif zq(s,x,y)\dif s\Bigg|\\
&\quad+\Bigg|\!\int^{\frac{t}{2}}_0\!\!\!\int_{\mR^d}\Delta^{\frac{\alpha+\gamma}{2}}p_0(t-s,x,z)q(s,z,y)\dif z\dif s\Bigg|\\
&=:\sC_1(t,x,y)+\sC_2(t,x,y)+\sC_3(t,x,y).
\end{align*}
For $\gamma<\beta$, we choose a $\gamma'>0$ such that $\gamma+\gamma'<\beta$, and by (\ref{e2}), (\ref{q2}) and (\ref{3p}), we have
\begin{align*}
\sC_1(t,x,y)&\leq C_1\!\!\int_{\frac{t}{2}}^t\!\!\!\int_{\mR^d}\varrho^{\beta-\gamma'}_{-\gamma}(t-s,x-z)\Big(\varrho^0_{\gamma'}+\varrho^{\beta}_{\gamma'-\beta}\Big)(s,z-y)\dif z\dif s\\
&\quad+C_1\!\!\int_{\frac{t}{2}}^t\!\!\!\int_{\mR^d}\varrho^{\beta-\gamma'}_{-\gamma}(t-s,x-z)\dif z\Big(\varrho^0_{\gamma'}+\varrho^{\beta}_{\gamma'-\beta}\Big)(s,x-y)\dif s\\
&\leq C_1\!\!\int_{0}^t\!\!\!\int_{\mR^d}\varrho^{\beta-\gamma'}_{-\gamma}(t-s,x-z)\Big(\varrho^0_{\gamma'}+\varrho^{\beta}_{\gamma'-\beta}\Big)(s,z-y)\dif z\dif s\\
&\quad+C_2\!\!\int_{\frac{t}{2}}^t(t-s)^{\frac{\beta-\gamma-\gamma'}{\alpha}-1}\Big(\varrho^0_{\gamma'}+\varrho^{\beta}_{\gamma'-\beta}\Big)(s,x-y)\dif s\\
&\leq C_3\Big(\varrho^{0}_{\beta-\gamma}+\varrho^{\beta-\gamma'}_{\gamma'-\gamma}+\varrho^{\beta}_{-\gamma}\Big)(t,x-y)\\
&\quad+C_3\!\!\int_{\frac{t}{2}}^t(t-s)^{\frac{\beta-\gamma-\gamma'}{\alpha}-1}\Big(\varrho^0_{\gamma'}+\varrho^{\beta}_{\gamma'-\beta}\Big)(s,x-y)\dif s\leq C_4\varrho^0_{-\gamma}(t,x-y).
\end{align*}
Thanks to (\ref{00}) and taken into account of (\ref{q1}), it holds
\begin{align*}
\sC_2(t,x,y)\leq C_1\!\!\int^t_{\frac{t}{2}}(t-s)^{\frac{\beta-\gamma}{\alpha}-1}\Big(\varrho_0^{\beta}+\varrho_\beta^0\Big)(s,x-y)\dif s\leq C_2\varrho^0_{-\gamma}(t,x-y).
\end{align*}
Finally, we have by (\ref{e2}), (\ref{q1}) and (\ref{3p}) that for any $\gamma{'\!'}>0$,
\begin{align*}
\sC_3(t,x,y)&\leq C_1\!\!\int^{\frac{t}{2}}_0\!\!\!\int_{\mR^d}\varrho^0_{-\gamma}(t-s,x-z)\Big(\varrho_0^{\beta}+\varrho_\beta^0\Big)(s,z-y)\dif z\dif s\\
&\leq C_2t^{-\frac{\gamma+\gamma{'\!'}}{\alpha}}\!\int^{t}_0\!\!\!\int_{\mR^d}\varrho^0_{\gamma{'\!'}}(t-s,x-z)\Big(\varrho_0^{\beta}+\varrho_\beta^0\Big)(s,z-y)\dif z\dif s\leq C_3\varrho^0_{-\gamma}(t,x-y).
\end{align*}
Based on the above estimates, we thus get (\ref{es}) by (\ref{e2}) and (\ref{heat}). Next, we proceed to prove (\ref{ph}). As in the proof of Lemma \ref{as}, we set
$$
\sK_{\nabla p_0}(t;x,x';y):=\nabla p_0(t,x,y)-\nabla p_0(t,x',y).
$$
Then, estimate (\ref{p0h}) yields that for any $\vartheta\in(0,1)$,
$$
|\sK_{\nabla p_0}(t;x,x';y)|\leq C_\vartheta|x-x'|^\vartheta\varrho_{\alpha-1-\vartheta}^0(t,\tilde x-y).
$$
As above, we can write
\begin{align*}
\big|\nabla\sS(t,x,y)-\nabla\sS(t,x',y)\big|&\leq \Bigg|\int_{\frac{t}{2}}^t\!\!\!\int_{\mR^d}\sK_{\nabla p_0}(t-s;x,x';z)\Big(q(s,z,y)-q(s,\tilde x,y)\Big)\dif z\dif s\Bigg|\\
&\quad+\Bigg|\int_{\frac{t}{2}}^t\!\!\!\int_{\mR^d}\sK_{\nabla p_0}(t-s;x,x';z)\dif zq(s,\tilde x,y)\dif s\Bigg|\\
&\quad+\Bigg|\!\int^{\frac{t}{2}}_0\!\!\!\int_{\mR^d}\sK_{\nabla p_0}(t-s;x,x';z)q(s,z,y)\dif z\dif s\Bigg|\\
&=:\sD_1(t,x,x',y)+\sD_2(t,x,x',y)+\sD_3(t,x,x',y).
\end{align*}
For $\vartheta<\alpha+\beta-1$, choose a $\gamma'$ such that $\vartheta+\gamma'<\alpha+\beta-1$, we follow the same procedure as in the estimate of $\sC_1(t,x,y)$ to derive that
\begin{align*}
\sD_1(t,x,x',y)&\leq C_\vartheta|x-x'|^\vartheta\!\!\int_{\frac{t}{2}}^t\!\!\!\int_{\mR^d}\varrho^{\beta-\gamma'}_{\alpha-1-\vartheta}(t-s,\tilde x-z)\Big(\varrho^0_{\gamma'}+\varrho^{\beta}_{\gamma'-\beta}\Big)(s,z-y)\dif z\dif s\\
&\quad+C_\vartheta|x-x'|^\vartheta\!\!\int_{\frac{t}{2}}^t\!\!\!\int_{\mR^d}\varrho^{\beta-\gamma'}_{\alpha-1-\vartheta}(t-s,\tilde x-z)\dif z\Big(\varrho^0_{\gamma'}+\varrho^{\beta}_{\gamma'-\beta}\Big)(s,x-y)\dif s\\
&\leq C_\vartheta|x-x'|^\vartheta\varrho_{\alpha-1-\vartheta}^0(t,\tilde x-y).
\end{align*}
As for $\sD_2(t,x,x',y)$ and $\sD_3(t,x,x',y)$, we may use (\ref{000}) and (\ref{p0h}) respectively, and argue the same way as estimating $\sC_2(t,x,y)$ and $\sC_3(t,x,y)$ to get that
$$
\sD_2(t,x,x',y)+\sD_3(t,x,x',y)\leq C_\vartheta|x-x'|^\vartheta\varrho_{\alpha-1-\vartheta}^0(t,\tilde x-y),
$$
which in turn yields (\ref{ph}).
\end{proof}

\section{Smoothing properties of the semigroup}
Let $\cT_t$ be the semigroup corresponding to $\sL^\kappa_{\alpha}$, that is,
$$
\cT_tf(x):=\int_{\mR^d}p(t,x,y)f(y)\dif y,\quad\forall f\in\cB_b(\mR^d).
$$
We shall use heat kernel estimates obtained in the last section to derive some space regularities of $\cT_t$, which has its own independent interest and will be used to study Krylov-type estimate and Zvonkin's transformation in the next section.

Let us first introduce some notations. Let $p\geq1$ and $\|\cdot\|_p$ denote the norm in $L^p(\mR^d)$. For $0<\gamma<2$, define the Bessel potential space $\mH^{\gamma}_p:=\mH^{\gamma}_p(\mR^d)$ by
$$
\mH^{\gamma}_p(\mR^d):=\big\{f\in L^p(\mR^d): \Delta^{\frac{\gamma}{2}}f\in L^p(\mR^d)\big\}
$$
with norm
$$
\|f\|_{\gamma,p}:=\|f\|_p+\|\Delta^{\frac{\gamma}{2}}f\|_p.
$$
In fact, this space can also be defined to be the complete space of $C^{\infty}_0(\mR^d)$ under the norm
$$
\|\sF^{-1}\big((1+|\cdot|^{\gamma})(\sF f)\big)\|_p<\infty,\quad \forall f\in C^{\infty}_0(\mR^d),
$$
where $\sF$ (resp. $\sF^{-1}$) denotes the Fourier transform (resp. the Fourier inverse transform). By Sobolev's embedding theorem, if $\gamma-\frac{d}{p}>0$ is not an integer, then (\cite[p. 206, (16)]{Tri})
\begin{align}
\mH^{\gamma}_p\hookrightarrow C^{\gamma-\frac{d}{p}}_b(\mR^d),\label{emb}
\end{align}
where for some $\gamma>0$, $C^{\gamma}_b(\mR^d)$ is the usual H\"older space with norm
$$
\|f\|_{C^{\gamma}_b}:=\sum_{i=1}^{[\gamma]}\|\nabla^i f(x)\|_\infty+\big[\nabla^{[\gamma]}f\big]_{\gamma-[\gamma]},
$$
here, $[\gamma]$ denotes the integer part of $\gamma$, and for a function $f$ on $\mR^d$ and $\vartheta\in(0,1)$,
\begin{align}
[f]_{\vartheta}:=\sup_{x\neq y}\frac{|f(x)-f(y)|}{|x-y|^{\vartheta}}. \label{hol}
\end{align}
Noticing that for $n\in \mN$, $\mH^{n}_p$ is just the usual Sobolev space with equivalent norm (\cite[p. 135, Theorem 3]{St})
$$
\|f\|_{n,p}=\|f\|_p+\|\nabla^nf\|_{p},
$$
here and below, $\nabla$ denotes the weak derivative of $f$. While for $0<\gamma\neq$ integer, the fractional Sobolev space $\mW^{\gamma}_p$ is defined by
\begin{align}
\|f\|_{\mW^{\gamma}_p}:=\|f\|_p+\sum_{k=0}^{[\gamma]}\Bigg(\int_{\mR^d}\!\!\int_{\mR^d}\frac{|\nabla^kf(x)-\nabla^kf(y)|^p}{|x-y|^{d+(\gamma-[\gamma])p}}\dif x\dif y\Bigg)^{1/p}<\infty.\label{fs}
\end{align}
The relation between $\mH^{\gamma}_p$ and $\mW^{\gamma}_p$ is that (cf. \cite[p. 190]{Tri}): for $\gamma>0$, $\eps\in(0,\gamma)$ and $p>1$,
\begin{align}
\mH^{\gamma+\eps}_p\hookrightarrow\mW^{\gamma}_p\hookrightarrow\mH^{\gamma-\eps}_p.\label{re}
\end{align}
Moreover, the following relationship can be found in \cite[p. 185]{Tri}: for $p>1$,  $\gamma_1\neq\gamma_2$ and $\vartheta\in(0,1)$,
\begin{align}
[\mH^{\gamma_1}_p,\mH^{\gamma_2}_p]_\vartheta=\mH^{\gamma_1+\vartheta(\gamma_2-\gamma_1)}_p,   \label{fw}
\end{align}
where $[A, B]_{\vartheta}$ denotes the complex interpolation space between two Banach space $A$ and $B$.
Recall the following complex interpolation result (\cite[p. 59, Theorem (a)]{Tri}).
\bl\label{inter}
Let $A_i\subseteq B_i$, $i=0,1$ be Banach spaces and $\sT: A_i\rightarrow B_i$, $i=0,1$ be a bounded linear operator. For any $\theta\in(0,1)$, we have
$$
\|\sT\|_{A_{\theta}\rightarrow B_{\theta}}\leq \|\sT\|_{A_0\rightarrow B_0}^{1-\theta}\|\sT\|_{A_1\rightarrow B_1}^{\theta},
$$
where $A_{\theta}:=[A_0, A_1]_{\theta}$, $B_{\theta}:=[B_0, B_1]_{\theta}$, and $\|\sT\|_{A_{\theta}\rightarrow B_{\theta}}$ denotes the operator norm of $\sT$ mapping $A_{\theta}$ to $B_{\theta}$.
\el

Given a locally integrable function $f$ on $\mR^d$, the Hardy-Littlewood maximal function of $f$ is defined by
$$
\cM f(x):=\sup_{0<r<\infty}\frac{1}{|B_r|}\int_{B_r}|f(x+y)|\dif y,
$$
where $|B_r|$ denotes the Lebesgue measure of ball $B_r$. The following well known result can be found in \cite[p. 5, Theorem 1]{St} and \cite{Zh3}.

\bl
(i) For $p\in(1,\infty]$ and all $f\in L^p(\mR^d)$, there exists a constant $C_{d,p}>0$ such that
\begin{align}
\|\cM f\|_p\leq C_{d,p}\|f\|_p.   \label{mf}
\end{align}
(ii) For every $f\in \mH^1_p$, there is a constant $C_d>0$ such that for a.e. $x,y\in\mR^d$,
\begin{align}
|f(x)-f(y)|\leq C_{d}|x-y|\Big(\cM|\nabla f|(x)+\cM|\nabla f|(y)\Big).   \label{w11}
\end{align}
\el

\vspace{2mm}

Now, we proceed to study the regularities of the semigroup $\cT_t$.

\bl[{\bf Regularity in H\"older space}]
There exist constants $C_1, C_2$ such that for any $\vartheta\in(0,1)$,
\begin{align}
\|\nabla \cT_tf\|_{\infty}\leq C_{1}t^{\frac{\vartheta-1}{\alpha}}[f]_\vartheta ,\quad  \forall \ f \in C_0^\infty(\mR^d), \label{33}
\end{align}
and for $\vartheta\in(0,1)$, $\vartheta'\in(0,\alpha+\beta-1)$,
\begin{align}
[\nabla \cT_tf]_{\vartheta'}\leq C_{2}t^{\frac{\vartheta-\vartheta'-1}{\alpha}}[f]_\vartheta,\quad  \forall \ f \in C_0^\infty(\mR^d),\label{34}
\end{align}
where $[\cdot]$ is defined by (\ref{hol}).
\el
\begin{proof}
Using (\ref{11}), we find that
$$
\int_{\mR^d}\nabla p(t,x,y)\dif y=0.
$$
Thus, we can write
\begin{align}
\nabla \cT_tf(x)=\int_{\mR^d}\nabla p(t,x,y)\big[f(y)-f(x)\big]\dif y. \label{cha}
\end{align}
Thus, by (\ref{na}), it is easy to find that for any $\vartheta\in(0,1)$,
$$
\|\nabla \cT_tf\|_{\infty}\leq C_1t^{-\frac{1}{\alpha}} [f]_\vartheta\!\!\int_{\mR^d}\frac{|x-y|^\vartheta}{(|x-y|+t^{1/\alpha})^{d+\alpha}}\dif y\leq C_1t^{\frac{\vartheta-1}{\alpha}}[f]_\vartheta.
$$
To prove (\ref{34}), we write
\begin{align*}
\nabla \cT_tf(x)-\nabla \cT_tf(x')&=\int_{\mR^d}\Big(\nabla p(t,x,y)-\nabla p(t,x',y)\Big)\big(f(y)-f(\tilde x)\big)\dif y,
\end{align*}
where $\tilde x$ is the one of the two points $x$ and $x'$ which is nearer to $y$. Taking into account of (\ref{ph}) we arrive that for $0<\vartheta'<\alpha+\beta-1$,
\begin{align*}
\nabla \cT_tf(x)-\nabla \cT_tf(x')&\leq C_2|x-x'|^{\vartheta'}[f]_\vartheta\!\int_{\mR^d}\varrho_{\alpha-1-\vartheta'}^\vartheta(t,\tilde x-y)\dif y\no\\
&\leq C_2|x-x'|^{\vartheta'}t^{\frac{\vartheta-\vartheta'-1}{\alpha}}[f]_\vartheta,
\end{align*}
which in turn implies the desired result.
\end{proof}

\bl[{\bf Regularity in Bessel potential space}]
Let $\theta\in(0,1)$ and $\gamma+\theta<\alpha+\beta$ hold, then for every $p>1$,
\begin{align}
\|\cT_tf\|_{\gamma+\theta,p}\leq C_{\gamma,p}t^{-\gamma/\alpha}\|f\|_{\theta,p}, \ \ \ \ \forall \ f \in \mH^{\theta}_p,    \label{tt}
\end{align}
where $C_{\gamma,p}>0$ is a constant.
\el
\begin{proof}
Thanks to a standard approximation argument, we only need to prove the estimate for $f\in C^{\infty}_0(\mR^d)$.
Observe
\begin{align*}
\|T_tf\|^p_{p}&=\int_{\mR^d}\Bigg(\int_{\mR^d}p(t,x,y)f(y)\dif y\Bigg)^p\dif x\\
&=\int_{\mR^d}\Bigg(\int_{\mR^d}p^{1/q}(t,x,y)p^{1/p}(t,x,y)f(y)\dif y\Bigg)^p\dif x\\
&\leq \int_{\mR^d}\Bigg(\int_{\mR^d}p(t,x,y)\dif y\Bigg)^{p/q}\Bigg(\int_{\mR^d}p(t,x,y)f^p(y)\dif y\Bigg)\dif x,
\end{align*}
thus
$$
\|\cT_tf\|_{p}\leq \|f\|_p.
$$
Hence, $\cT_t$ is a contraction semigroup and we may assume $t<1$ below. By Fubini's theorem and \eqref{es}, for $\beta_1<\alpha+\beta$ we have
\begin{align*}
\Delta^{\frac{\beta_{1}}{2}} \cT_tf(x)&=\int_{\mR^d}\Delta^{\frac{\beta_{1}}{2}}p(t,x,y)f(y)\dif y\\
&\leq C_1t^{-\beta_1/\alpha}\int_{\mR^d}\varrho_{\alpha}^0(t,x-y)f(y)\dif y,
\end{align*}
which implies that for $\beta_1<\alpha+\beta$,
$$
\|\cT_tf\|_{\beta_1,p}\leq C_1t^{-\beta_1/\alpha}\|f\|_p.
$$
Let $\beta_2$ be such that $1+\beta_2<\alpha+\beta$, as in (\ref{cha}), we write
\begin{align*}
\Delta^{\frac{1+\beta_2}{2}} \cT_tf(x)&=\int_{\mR^d}\Delta^{\frac{1+\beta_2}{2}}p(t,x,y)\big[f(y)-f(x)\big]\dif y.
\end{align*}
Then, it follows by (\ref{es}) and (\ref{w11}) that
\begin{align*}
|\Delta^{\frac{1+\beta_2}{2}} \cT_tf(x)|&\leq \int_{\mR^d}|\Delta^{\frac{1+\beta_2}{2}}p(t,x,y)|\cdot|x-y|\Big(\cM|\nabla f|(x)+\cM|\nabla f|(y)\Big)\dif y\\
&\leq C\!\int_{\mR^d}\varrho_{\alpha-1-\beta_2}^0(t,x-y)\cdot|x-y|\Big(\cM|\nabla f|(x)+\cM|\nabla f|(y)\Big)\dif y\\
&=C\cM|\nabla f|(x)\!\int_{\mR^d}\varrho_{\alpha-1-\beta_2}^0(t,y)|y|\dif y+C\!\int_{\mR^d}\varrho_{\alpha-1-\beta_2}^0(t,y)|y|\cdot\cM|\nabla f|(x-y)\dif y.
\end{align*}
Since $\alpha>1$, one can check easily that
\begin{align*}
\int_{\mR^d}\varrho_{\alpha-1-\beta_2}^0(t,y)|y|\dif y=
\left(\int_{|y| \le t^{1/\alpha}}+\int_{|y|>t^{1/\alpha}}\right)\varrho_{\alpha-1-\beta_2}^0(t,y)|y|\dif y\leq t^{-\frac{\beta_2}{\alpha}},
\end{align*}
which, together with (\ref{mf}) and Minkovski's inequality, yields
\begin{align*}
\|\Delta^{\frac{1+\beta_2}{2}} \cT_tf\|_p&\leq Ct^{-\frac{\beta_2}{\alpha}}\|\cM|\nabla f|\|_p+C\bigg\|\int_{\mR^d}\varrho_{\alpha-1-\beta_2}^0(t,y)|y|\cM|\nabla f|(x-y)\dif y\bigg\|_p\\
&\leq Ct^{-\frac{\beta_2}{\alpha}}\|\nabla f\|_p+C\!\!\int_{\mR^d}\varrho_{\alpha-1-\beta_2}^0(t,y)|y|\cdot\|\cM|\nabla f|(x-y)\|_p\dif y\\
&\leq Ct^{-\frac{\beta_2}{\alpha}}\|\nabla f\|_p.
\end{align*}
Hence, we get
$$
\|\cT_tf\|_{1+\beta_2,p}\leq Ct^{-\frac{\beta_2}{\alpha}}\|f\|_{1,p}.
$$
By the interpolation (\ref{fw}), for $\theta\in(0,1)$
$$
\mH^{\theta}_p=[L^p,\mH^{1}_p]_{\theta},\quad\mH^{\beta_1+(1+\beta_2-\beta_1)\theta}_p=[\mH^{\beta_1}_p,\mH^{1+\beta_2}_p]_{\theta},
$$
and Lemma \ref{inter}, we can derive that for $\gamma=(1-\theta)\beta_1+\theta\beta_2<\alpha+\beta-\theta$,
$$
\|\cT_tf\|_{\gamma+\theta,p}\leq C t^{-\frac{\gamma}{\alpha}}\|f\|_{\theta,p}.
$$
The proof is finished.
\end{proof}

\br\label{ree}
The proof in \cite{Pri} and \cite{Zh00} relies heavily on the symmetry of $\Delta^{\frac{\alpha}{2}}$ and the smooth properties at least up to second order of its heat kernel $p_{\alpha}(t,x)$. However, the operator $\sL_{\alpha}^\kappa$ considered here is non-symmetric and its heat kernel has no more regularity than `$\alpha+\beta$'-order, as we have seen in Lemma \ref{seen}.
\er

\section{Krylov-type estimate and Zvonkin's transformation}
This section consists of two parts, one is to obtain a Krylov-type estimate for the strong solution of SDE \eqref{sde2}, while the other is to transforms SDE (\ref{sde2}) into a new one with better drift coefficient by Zvonkin's transformation. The regularity of the semigroup $\cT_t$ obtained in the last section will play an important role in the two subsections below.

\subsection{Krylov's estimate}

Let $X_t(x)$ be a strong solution to SDE (\ref{sde2}). Usually, the It\^o's formula is performed for functions $f\in C^2_b(\mR^d)$. However, this is too strong for our latter use. Notice that by (\ref{s1}), $\sL_{\alpha}^\kappa f$ is meaningful for any $f\in C^{\gamma}_b(\mR^d)$ as long as $\gamma>\alpha$. Indeed, we have by (\ref{nu}) that
\begin{align*}
\sL^\kappa_{\alpha}f(x)&\leq C_{d,\alpha}\!\int_{|z|\leq 1}\!\!\int_0^1\!\big|\nabla f(x+rz)-\nabla f(x)\big|\dif r\frac{\dif z}{|z|^{d+\alpha-1}}+C_{d,\alpha}\|f\|_{\infty}\\
&\leq C_{d,\alpha}\!\int_{|z|\leq 1}\frac{\dif z}{|z|^{d+\alpha-\gamma}}\|f\|_{\gamma}+C_{d,\alpha}\|f\|_{\infty}<\infty.
\end{align*}
We first show that It\^o's formula holds for $f(X_t)$ when $f\in C^{\gamma}_b(\mR^d)$ with $\gamma>\alpha$.
\bl\label{ito}
Suppose (\ref{nu}) and (\ref{s1}) hold. Let $X_t$ satisfy (\ref{sde2}) and $f\in C^{\gamma}_b(\mR^d)$ with $\gamma>\alpha$. Then, we have
\begin{align*}
f(X_t)-f(x)-\int_0^t\!\sL f(X_s)\dif s=\int_0^t\!\!\!\int_0^{\infty}\!\!\!\!\int_{\mR^d}\big[f\big(X_{s-}+1_{[0,\sigma(X_{s-},z)]}(r)z\big)-f(X_s)\big]\tilde \cN(\dif z\times\dif r\times \dif s),
\end{align*}
where $\sL=\sL_\alpha^\kappa+b\cdot\nabla$.
\el
\begin{proof}
Let $\rho\in C^{\infty}_0(\mR^d)$ such that $\int_{\mR^d}\rho(x)\dif x=1$. Define $\rho_n(x):=n^d\rho(nx)$, and
\begin{align}
f_n(x):=\int_{\mR^d}f(y)\rho_n(x-y)\dif y.\label{mo}
\end{align}
Hence, we have $f_n\in C^2_b(\mR^d)$ with $\|f_n\|_{C^\gamma_b}\leq \|f\|_{C^\gamma_b}$, and $\|f_n-f\|_{C^{\gamma'}_b}\rightarrow0$ for every $\gamma'<\gamma$. By using It\^o's formula for $f_n(X_t)$, we get
\begin{align*}
f_n(X_t)-f_n(x)-\int_0^t\!\sL f_n(X_s)\dif s=\int_0^t\!\!\!\int_0^{\infty}\!\!\!\!\int_{\mR^d}\big[f_n\big(X_{s-}+1_{[0,\sigma(X_{s-},z)]}(r)z\big)-f_n(X_{s-})\big]\tilde \cN(\dif z\times\dif r\times \dif s).
\end{align*}
Now we are going to pass the limits on the both sides of the above equality. It is easy to see that for every $\omega$ and $x\in\mR^d$,
$$
f_n(X_t)-f_n(x)\rightarrow f(X_t)-f(x),\quad \text{as}\,\,n\rightarrow\infty.
$$
Since
\begin{align*}
|f_n(x+z)-f_n(x)-z\cdot\nabla f_n(x)|\leq C|z|^{\gamma}\|f_n\|_{C^\gamma_b}\leq C|z|^{\gamma}\|f\|_{C^\gamma_b},
\end{align*}
we can get by dominated convergence theorem that for every $\omega$,
$$
\int_0^t\!\sL f_n(X_s)\dif s\rightarrow\int_0^t\!\sL f(X_s)\dif s,\quad \text{as}\,\,n\rightarrow\infty.
$$
Finally, by the isometry formula, we have
\begin{align*}
&\mE\bigg|\!\int_0^t\!\!\!\int_0^{\infty}\!\!\!\!\int_{\mR^d}\Big[f_n\big(X_{s-}+1_{[0,\sigma(X_{s-},z)]}(r)z\big)-f_n(X_s)\\ &\qquad-f\big(X_{s-}+1_{[0,\sigma(X_{s-},z)]}(r)z\big)+f(X_{s-})\Big]\tilde \cN(\dif z\times\dif r\times \dif s)\bigg|^2\\
&=\mE\int_0^t\!\!\!\int_{\mR^d}\!\!\int_0^{\infty}1_{[0,\sigma(X_s,z)]}(r)\big|f_n(X_s+z)-f_n(X_s) -f(X_s+z)+f(X_s)\big|^2\dif r\nu(\dif z)\dif s\\
&\leq C\!\!\int_0^t\!\!\!\int_{\mR^d}\mE\big|f_n(X_s+z)-f_n(X_s) -f(X_s+z)+f(X_s)\big|^2\nu(\dif z)\dif s\rightarrow0,\quad \text{as}\,\,n\rightarrow\infty,
\end{align*}
where in the last step we have used the fact that $\sigma$ is bounded, $\|f_n\|_{C^\gamma_b}\leq \|f\|_{C^\gamma_b}$ and the dominated convergence theorem again. The proof is finished.
\end{proof}

We need the following results about the semi-linear elliptic PDE to prove the Krylov's estimate.

\bl\label{p1}
Let $\lambda,\bar k>0$. For any $f\in C_b^\infty(\mR^d)$, there exists a unique classical solution $u\in C_b^{\alpha+\theta}(\mR^d)$ with $\theta<\beta$ to equation
\begin{align}
\lambda u-\sL^\kappa_\alpha u-\bar k|\nabla u|=f,   \label{pide3}
\end{align}
which also satisfies the following integral equation:
\begin{align}
u(x)=\int_0^\infty\!\e^{-\lambda t}\cT_t\big(\bar k|\nabla u|+f\big)(x)\dif t.\label{u}
\end{align}
Moreover, for $\lambda$ big enough, we have for any $1<\gamma<\alpha$,
\begin{align}
\|u\|_{\gamma,p}\leq C\|f\|_p.\label{ess}
\end{align}
\el
\begin{proof}
Let us first construct the solution of \eqref{u} via Picard's iteration argument.
Set $u_0\equiv0$ and for $n\in \mN$, define $u_n$ recursively by
$$
u_n(x):=\int_0^\infty\!\e^{-\lambda t}\cT_t\big(\bar k|\nabla u_{n-1}|+f\big)(x)\dif t.
$$
In view of (\ref{na}), it is easy to check that $u_1\in C^1_b(\mR^d)$, and $u_2$ is thus well defined, and so on. We further write that
\begin{align*}
\nabla u_1(x)-\nabla u_1(y)&=\int_0^\infty\!\e^{-\lambda t}\big[\nabla \cT_tf(x)-\nabla \cT_tf(y)\big]\dif t\\
&=\left(\int_0^{|x-y|^{\alpha}}+\int^{\infty}_{|x-y|^{\alpha}}\right)\e^{-\lambda t}\big[\nabla \cT_tf(x)-\nabla \cT_tf(y)\big]\dif t=:I_1+I_2.
\end{align*}
For $I_1$, we have by (\ref{33}) that for $0<\vartheta<2-\alpha$,
\begin{align*}
I_1\leq C\!\!\int_0^{|x-y|^{\alpha}}\!t^{\frac{\vartheta-1}{\alpha}}\dif t \cdot[f]_\vartheta\leq C|x-y|^{\alpha+\vartheta-1}[f]_\vartheta.
\end{align*}
Meanwhile, as a direct result of (\ref{34}), we have for $\vartheta\in(0,\beta)$ and $\alpha+\vartheta-1<\vartheta'<\alpha+\beta-1$,
\begin{align*}
I_2\leq C|x-y|^{\vartheta'} [f]_\vartheta\!\int^{\infty}_{|x-y|^{\alpha}}t^{\frac{\vartheta-\vartheta'-1}{\alpha}}\dif t\leq C|x-y|^{\alpha+\vartheta-1} [f]_\vartheta.
\end{align*}
Consequently, we get that $u_1\in C_b^{\alpha+\vartheta}(\mR^d)$ if $f\in C_b^{\vartheta}(\mR^d)$. Noticing that $u_1\in C_b^{\alpha+\vartheta}(\mR^d)$ implies that
$|\nabla u_1|\in C_b^{\vartheta}(\mR^d)$ since
$$
\big||\nabla u_1|(x)-|\nabla u_1|(y)\big|\leq |\nabla u_1(x)-\nabla u_1(y)|.
$$
Repeating the above argument, we have for every $n\in\mN$ and $\vartheta\in(0,\beta)$,
$$
u_n\in C_b^{\alpha+\vartheta}(\mR^d).
$$
Moreover, since
\begin{align*}
u_n(x)-u_m(x)&=\int_0^\infty\!\e^{-\lambda t}\cT_t\big(\bar k|\nabla u_{n-1}|-k|\nabla u_{m-1}|\big)(x)\dif t\\
&\leq \bar k\!\!\int_0^\infty\!\e^{-\lambda t}\cT_t\big|\nabla u_{n-1}-\nabla u_{m-1}\big|(x)\dif t,
\end{align*}
we further have that for $\vartheta'\!'$ with $\vartheta<\vartheta'\!'<\alpha+\vartheta-1$,
\begin{align*}
\|u_n-u_m\|_{C_b^{\alpha+\vartheta}}&\leq \bar k\!\!\int_0^\infty\!\e^{-\lambda t}\|\cT_t\big(|\nabla u_{n-1}-\nabla u_{m-1}|\big)\|_{C_b^{\alpha+\vartheta}}\dif t\\
&\leq C_{\bar k}\!\!\int_0^\infty\!\e^{-\lambda t}\Big(1+t^{\frac{\vartheta-1}{\alpha}}+t^{\frac{\vartheta'\!'-\vartheta}{\alpha}-1}\Big)\dif t\cdot\|u_{n-1}-u_{m-1}\|_{C_b^{1+\vartheta'\!'}}\\
&\leq C_{\bar k}\lambda^{-\frac{\vartheta'\!'-\vartheta}{\alpha}}\|u_{n-1}-u_{m-1}\|_{C_b^{\alpha+\vartheta}},
\end{align*}
where we have also used the fact that $\lambda>1$. This means that $u_n$ is Cauchy sequence in $C_b^{\alpha+\vartheta}(\mR^d)$. Thus, there exists a limit function $u\in C_b^{\alpha+\vartheta}(\mR^d)$ with $\vartheta\in(0,\beta)$ satisfying (\ref{u}).

Since by \cite[(1.7)]{Ch-Zh}, for every $f\in C^{\vartheta}_b(\mR^d)$,
$$
\p t\cT_tf(x)=\sL^\kappa_{\alpha}\cT_tf(x),
$$
we have by integral by part,
\begin{align*}
\sL^\kappa_{\alpha}u(x)&=\int_0^\infty\!\e^{-\lambda t}\sL^\kappa_{\alpha}\cT_t\big(\bar k|\nabla u|+f\big)(x)\dif t\\
&=\int_0^\infty\!\e^{-\lambda t}\p t\cT_t\big(\bar k|\nabla u|+f\big)(x)\dif t\\
&=-\bar k|\nabla u|(x)-f(x)+\lambda u(x),
\end{align*}
which means that $u$ satisfies PIDE (\ref{pide3}). Moreover, we have by (\ref{tt}) that
\begin{align*}
\|u\|_{\gamma,p}&\leq \int_0^\infty\!\e^{-\lambda t}\|\cT_t\big(\bar k|\nabla u|+f\big)\|_{\gamma,p}\dif t\\
&\leq C_{\bar k}\!\!\int_0^\infty\!\e^{-\lambda t}t^{-\frac{\gamma}{\alpha}}\dif t\Big(\|\nabla u\|_{p}+\|f\|_{p}\Big)\\
&\leq C_{\bar k}\lambda^{\frac{\gamma}{\alpha}-1}\Big(\|u\|_{1,p}+\|f\|_{p}\Big).
\end{align*}
Choosing $\lambda$ big enough such that $C_{\bar k}\lambda^{\frac{\gamma}{\alpha}-1}<1$, we can get (\ref{ess}). The whole proof is finished.
\end{proof}

Now, we can give the Krylov's estimate for strong solutions of SDE (\ref{sde2}).
\bl
Let $X_t$ be a strong solution of SDE (\ref{sde2}). Then, for any $T>0$, there exist a constant $C_T$ such that for any $f\in L^p(\mR^d)$ with $p>d/\alpha$, we have
\begin{align}
\mE\Bigg(\int_0^T\!\!f(X_s)\dif s\Bigg)\leq C_T \|f\|_p.  \label{kry1}
\end{align}
\el
\begin{proof}
We first suppose that $f\in C^{\infty}_0(\mR^{d})$. By Lemma \ref{p1}, there exists a solution $u\in C_b^{\alpha+\vartheta}(\mR^d)$ with $0<\vartheta<\beta$ to equation (\ref{pide3}), which is given by (\ref{u}). According to Lemma \ref{ito}, we can use It\^o's formula to get for any $t>0$,
\begin{align*}
u(X_t)&=u(x)+\int_0^t\!\!\sL^\kappa_{\alpha}u(X_s)+b\cdot\nabla u(X_s)\dif s\\ &\quad+\int_0^t\!\!\!\int_0^{\infty}\!\!\!\!\int_{\mR^d}\big[u\big(X_{s-}+1_{[0,\sigma(X_{s-},z)]}(r)z\big)-u(X_{s-})\big]\tilde \cN(\dif z\times\dif r\times \dif s).
\end{align*}
By (\ref{pide3}) and take $\bar k$ big enough such that $\bar k>\|b\|_{\infty}$, we have for any $T>0$,
\begin{align*}
\int_0^T\!\!\sL^\kappa_{\alpha}u(X_s)+b\cdot\nabla u(X_s)\dif s&=\lambda\int_0^T\!u(X_s)\dif s+\int_0^T\!\!\Big(b\cdot\nabla u(X_s)-\bar k|\nabla u|(X_s)\Big)\dif s-\int_0^T\!\!f(X_s)\dif s\\
&\leq \lambda\int_0^T\!u(X_s)\dif s-\int_0^T\!\!f(X_s)\dif s.
\end{align*}
Consequently, it follows that
\begin{align*}
\mE\Bigg(\int_0^T\!\!f(X_s)\dif s\Bigg)\leq \lambda\mE\Bigg(\int_0^T\!|u(X_s)|\dif s\Bigg)+2\|u\|_{\infty}\leq C_{\lambda,T}\|u\|_{\infty}.
\end{align*}
Since $p>d/\alpha$, there exists a $\gamma<\alpha$ such that $p>d/\gamma$. Now, using (\ref{ess}), we conclude that
\begin{align*}
\mE\Bigg(\int_0^T\!\!f(X_s)\dif s\Bigg)\leq C_{\lambda,T}\|u\|_{\infty}\leq C_{\lambda,T}\|u\|_{\gamma,p}\leq C_{\lambda,T}\|f\|_{p},
\end{align*}
where the second inequality is due to the Sobolev embedding (\ref{emb}).
By a standard density argument as in \cite{Zh00}, we get the desired result for general $f\in L^p(\mR^d)$.
\end{proof}

\subsection{Zvonkin's transformation}

Now, we follow the idea in \cite{Kr-Ro,Zh1,Zh00} to transform SDE (\ref{sde2}) into a new one
with better coefficients.  Unlike \cite{Kr-Ro,Zh1,Zh00}, we do not have  well developed elliptic equation theory to solve the following equation \eqref{pide1} in Bessel potential space (or Sobolev space):
\begin{align}
\lambda u(x)-\sL^\kappa_{\alpha} u(x)-b\cdot\nabla u(x)=b(x).   \label{pide1}
\end{align}
The point is that the operator $\sL^\kappa_{\alpha}$ is non-symmetric and has no comparison with operator $\Delta^{\frac{\alpha}{2}}$. To be more precisely, even if we know $\Delta^{\frac{\alpha}{2}}u\in L^p(\mR^d)$ for some $p>1$, it is not easy to claim that $\sL^\kappa_{\alpha}u\in L^p(\mR^d)$, and vice versa. The authors in \cite{Kr-Ro,Zh1} encountered with classical local seconder order differential operator
\begin{align*}
\sL_2:=\frac{1}{2}\sum_{i,j=1}^da_{ij}(x)\frac{\p^2}{\p{x_i}\p{x_j}},
\end{align*}
which has been well studied, it is clear that if $\big(a_{ij}(x)\big)$ is bounded, then $\Delta u\in L^p(\mR^d)$ implies $\sL_2 u\in L^p(\mR^d)$ for every $p>1$. While \cite{Zh00} only needs to handle the symmetric operator $\Delta^{\frac{\alpha}{2}}$. Therefore, some new ways are required in our case.

First, we note that the elliptic equation can still be solved in the framework of H\"older space. The following result can be proved as in Lemma \ref{p1}, we omit the details here.

\bl\label{int}
Assume that for some $\vartheta\in(0,\beta)$, $b\in C^\vartheta_b(\mR^d)$ . Then, there exists a classical solution $u\in C_b^{\alpha+\vartheta}(\mR^d)$ to (\ref{pide1}). Meanwhile, $u$ also satisfies the integral equation:
\begin{align}
u(x)=\int_0^\infty\!\e^{-\lambda t}\cT_t\big(b\cdot\nabla u+b\big)(x)\dif t.   \label{g1}
\end{align}
\el

Despite that we can not solve the elliptic equation (\ref{pide1}) in Bessel potential space (or Sobolev space), we can solve the integral equation (\ref{g1}) in this framework thanks to the Bessel regularity of $\cT_t$ obtained in Section 4.

\bt\label{equ}
Let $1<\gamma<\alpha$. Suppose that for some $p>\frac{d}{\gamma}$ and $0<\theta\in(1-\gamma+\frac{d}{p},1)$,
$$
b\in L^{\infty}(\mR^d)\cap\mW^{\theta}_p.
$$
Then, for $\lambda$ big enough there exists a function $u\in \mH^{\gamma+\theta}_p$ satisfying the integral equation (\ref{g1}).
Moreover, we have
\begin{align}
\|u\|_{\gamma+\theta,p}\leq C_1\|b\|_{\mW^\theta_p} \label{esa}
\end{align}
and
\begin{align}
\|u\|_{\infty}+\|\nabla u\|_{\infty}\leq \frac{1}{2}.    \label{es2}
\end{align}
\et
\begin{proof}
We only show the priori estimate (\ref{esa}). Then ,the existence of solutions follows by the standard continuity method. Since $1<\gamma<\alpha$, choose $\eps\in(0,\alpha-\gamma)$, we have by (\ref{tt}) and (\ref{re}) that
\begin{align*}
\|u\|_{\gamma+\theta,p}&\leq \int_0^\infty\!\e^{-\lambda t}\|\cT_{t}(b\cdot\nabla u+b)\|_{\gamma+\theta,p}\dif t\\
&\leq C_1\!\int_0^\infty\!\e^{-\lambda t}t^{-\frac{\gamma+\eps}{\alpha}}\dif t\cdot\Big(\|b\cdot\nabla u\|_{\theta-\eps,p}+\|b\|_{\theta-\eps,p}\Big)\\
&\leq C_1\lambda^{\frac{\gamma+\eps}{\alpha}-1}\cdot\Big(\|b\cdot\nabla u\|_{\mW^\theta_p}+\|b\|_{\mW^\theta_p}\Big).
\end{align*}
In view of (\ref{fs}), (\ref{emb}), and thanks to the condition that $\gamma>1$ and $\gamma+\theta-1-\tfrac{d}{p}>0$, we can know
\begin{align*}
\|b\cdot\nabla u\|_{\mW^\theta_p}&\leq \|b\|_\infty\|\nabla u\|_p+\|b\|_{\mW^\theta_p}\|\nabla u\|_{\infty}+\|b\|_{\infty}\|\nabla u\|_{\mW^\theta_p}\\
&\leq \|b\|_\infty\|u\|_{1,p}+\|b\|_{\mW^\theta_p}\|u\|_{\gamma+\theta,p}+ \|b\|_{\infty}\|u\|_{\mW^{1+\theta}_p}\\
&\leq \Big(\|b\|_{\mW^\theta_p}+\|b\|_{\infty}\Big)\|u\|_{\gamma+\theta,p},
\end{align*}
where we also used (\ref{re}) in the last inequality. It then follows that
$$
\|u\|_{\gamma+\theta,p}\leq C_1\lambda^{\frac{\gamma+\eps}{\alpha}-1}\Big(\|b\|_{\mW^\theta_p}+\|b\|_{\infty}\Big)\|u\|_{\gamma+\theta,p}+C_1\lambda^{\frac{\gamma+\eps}{\alpha}-1}\|b\|_{\mW^\theta_p}.
$$
Hence, we can choose $\lambda_1$ big enough such that
$$
C_1\lambda^{\frac{\gamma+\eps}{\alpha}-1}_1\Big(\|b\|_{\mW^\theta_p}+\|b\|_{\infty}\Big)<\frac{1}{2},
$$
which means (\ref{esa}) is true. Moreover, we can take $\lambda\geq \lambda_1$ such that
$$
\lambda^{\frac{\gamma+\eps}{\alpha}-1}\|b\|_{\mW^\theta_p}<\frac{1}{8},
$$
and then get
$$
\|u\|_{\infty}+\|\nabla u\|_{\infty}\leq2\|u\|_{\gamma+\theta,p}\leq\frac{1}{2}.
$$
The proof is finished.
\end{proof}

In the sequel, we assume that $b\in  L^{\infty}(\mR^d)\cap\mW^{\theta}_p$ with
\begin{align}
\theta\in(1-\frac{\alpha}{2},1),\quad p>\frac{2d}{\alpha}. \label{index}
\end{align}
Notice that we can always choose a $1<\gamma<\alpha$ such that
$$
\theta>1-\gamma+\frac{d}{p}\quad\text{and}\quad p>\frac{d}{\gamma}.
$$
Hence, according to Theorem \ref{equ}, for $\lambda$ big enough we can get a function $u\in \mH^{\gamma+\theta}_p$ satisfying the integral equation (\ref{g1}). Define
$$
\Phi(x):=x+u(x).
$$
In view of (\ref{es2}), we also have
$$
\frac{1}{2}|x-y|\leq\big|\Phi(x)-\Phi(y)\big|\leq \frac{3}{2}|x-y|,
$$
which implies that the map $x\rightarrow\Phi(x)$ forms a $C^1$-diffeomorphism and
\begin{align}
\frac{1}{2}\leq \|\nabla\Phi\|_{\infty},\|\nabla\Phi^{-1}\|_{\infty}\leq 2,   \label{upd}
\end{align}
where $\Phi^{-1}(\cdot)$ is the inverse function of $\Phi(\cdot)$.

We prove the following Zvonkin's transformation.

\bl\label{zvon}
Let $\Phi(x)$ be defined as above and $X_t$ solve SDE (\ref{sde2}). Then, $Y_t:=\Phi(X_t)$ satisfies the following SDE:
\begin{align}
Y_t&=\Phi(x)+\int_0^t\tilde b(Y_s)\dif s+\int_0^t\!\!\!\int_0^{\infty}\!\!\!\!\int_{|z|\leq 1}\tilde g(Y_{s-},z)1_{[0,\tilde\sigma(Y_{s-},z)]}(r)\tilde \cN(\dif z\times\dif r\times \dif s)\no\\
&\quad+\int_0^t\!\!\!\int_0^{\infty}\!\!\!\!\int_{|z|> 1}\tilde g(Y_{s-},z)1_{[0,\tilde\sigma(Y_{s-},z)]}(r)\cN(\dif z\times\dif r\times \dif s), \label{sde3}
\end{align}
where
\begin{align*}
\tilde b(x)=\lambda u\big(\Phi^{-1}(x)\big)-\int_{|z|>1}\!\big[u\big(\Phi^{-1}(x)+z\big)-u\big(\Phi^{-1}(x)\big)\big]\sigma\big(\Phi^{-1}(x),z\big)\nu(\dif z)
\end{align*}
and
\begin{align*}
\tilde g(x,z):=\Phi\big(\Phi^{-1}(x)+z\big)-x,\quad \tilde\sigma(x,z):=\sigma\big(\Phi^{-1}(x),z\big).
\end{align*}
\el
\begin{proof}
Let $b_n$ be the mollifying approximation of $b$ defined as in (\ref{mo}). Then, it is obvious that
\begin{align}
\|b_n\|_{\infty}\leq \|b\|_{\infty},\quad \|b_n-b\|_{\mW^{\vartheta}_p}\rightarrow0,\,\,\text{as}\,\,n\rightarrow\infty. \label{bn}
\end{align}
Meanwhile, $b_n\in C_b^\vartheta(\mR^d)$ for any $\vartheta\in(0,\beta)$ and we may assume
\begin{align}
b_n\rightarrow b,\,\,a.e.,\,\,\text{as}\,\,n\rightarrow\infty. \label{bnn}
\end{align}
Let $u_n\in C_b^{\alpha+\vartheta}(\mR^d)$ be the classical solution to the elliptic equation (\ref{pide1}) with $b$ replaced by $b_n$. According to Lemma \ref{int}, we know that $u_n$ also satisfies the integral equation
\begin{align*}
u_n(x)=\int_0^\infty\!\e^{-\lambda t}\cT_t\big(b_n\cdot\nabla u_n+b_n\big)(x)\dif t.
\end{align*}
We proceed to show that
\begin{align}
\|u_n-u\|_{\gamma+\vartheta,p}\rightarrow0,\,\,\text{as}\,\,n\rightarrow\infty. \label{unn}
\end{align}
In fact, write
\begin{align*}
u_n(x)-u(x)=\int_0^\infty\!\e^{-\lambda t}\cT_t\Big(b_n(\nabla u_n-\nabla u)+(b_n-b)\cdot\nabla u+(b_n-b)\Big)(x)\dif t.
\end{align*}
As in the proof of Theorem \ref{equ}, we have
\begin{align*}
\|u_n-u\|_{\gamma+\theta,p}\leq C_1\lambda^{\frac{\gamma+\eps}{\alpha}-1}\cdot\Big(\|b_n(\nabla u_n-\nabla u)\|_{\mW^\theta_p}+\|(b_n-b)\cdot\nabla u\|_{\mW^\theta_p}+\|b_n-b\|_{\mW^\theta_p}\Big).
\end{align*}
At the same time, we can also get by (\ref{bn}) that
\begin{align*}
\|b_n(\nabla u_n-\nabla u)\|_{\mW^\theta_p}\leq \Big(\|b_n\|_{\mW^\theta_p}+\|b_n\|_{\infty}\Big)\|u_n-u\|_{\gamma+\theta,p}\leq \Big(\|b\|_{\mW^\theta_p}+\|b\|_{\infty}\Big)\|u_n-u\|_{\gamma+\theta,p}.
\end{align*}
Hence,
\begin{align*}
\|u_n-u\|_{\gamma+\theta,p}&\leq C_\lambda\Big(\|(b_n-b)\cdot\nabla u\|_{\mW^\theta_p}+\|b_n-b\|_{\mW^\theta_p}\Big)\\
&\leq C_\lambda\|b_n-b\|_{\mW^\theta_p}+C_\lambda\!\!\int_{\mR^d}\!\!\int_{\mR^d}\frac{|\nabla u(x)-\nabla u(y)|^p}{|x-y|^{d+\vartheta p}}|b_n(y)-b(y)|\dif x\dif y\rightarrow0,\,\,\text{as}\,\,n\rightarrow\infty,
\end{align*}
where we have used (\ref{bn}), (\ref{bnn}) and the dominated convergence theorem.

Now, we define
$$
\Phi_n(x):=x+u_n(x).
$$
By Lemma \ref{ito} and recalling that $u_n$ satisfies (\ref{pide1}), we can use It\^o's formula for $u_n$ to get
\begin{align*}
u_n(X_t)&=u_n(x)+\int_0^t\!\Big(\lambda u_n+\big(b-b_n\big)\cdot\nabla u_n-b_n\Big)(X_s)\dif s\\
&\quad-\int_0^t\!\!\!\int_0^{\infty}\!\!\!\int_{|z|>1}\big[u_n\big(X_s+1_{[0,\sigma(X_s,z)]}(v)z\big)-u_n(X_s)\big]\nu(\dif z)\dif r\dif s\\
&\quad+\int_0^t\!\!\!\int_0^{\infty}\!\!\!\!\int_{|z|> 1}\big[u_n\big(X_{s-}+1_{[0,\sigma(X_{s-},z)]}(v)z\big)-u_n(X_{s-})\big] \cN(\dif z\times\dif r\times \dif s)\\
&\quad+\int_0^t\!\!\!\int_0^{\infty}\!\!\!\!\int_{|z|\leq 1}\big[u_n\big(X_{s-}+1_{[0,\sigma(X_{s-},z)]}(v)z\big)-u_n(X_{s-})\big]\tilde \cN(\dif z\times\dif r\times \dif s).
\end{align*}
Adding this with SDE \eqref{sde2} and noticing that
$$
\Phi_n\big(x+y\big)-\Phi_n(x)=u_n(x+y)-u_n(x)+y,
$$
\begin{align*}
\Phi_n\big(X_{s-}+1_{[0,\sigma(X_{s-},z)]}(r)z\big)-\Phi_n(X_{s-})=1_{[0,\sigma(X_{s-},z)]}(r)\left[\Phi_n\big(X_{s-}+z\big)-\Phi_n(X_{s-})\right],\
\end{align*}
we obtain
\begin{align*}
Y^n_t&:=\Phi_n(X_t)=\Phi_n(x)+\int_0^t\!\lambda u_n(X_s)\dif s+\int_0^t\!\Big(\big(b-b_n\big)\cdot\nabla u_n+\big(b-b_n\big)\Big)(X_s)\dif s\\
&\quad-\int_0^t\!\!\!\int_0^{\infty}\!\!\!\int_{|z|>1}\big[u_n\big(X_s+1_{[0,\sigma(X_s,z)]}(r)z\big)-u_n(X_s)\big]\nu(\dif z)\dif r\dif s\\
&\quad+\int_0^t\!\!\!\int_0^{\infty}\!\!\!\!\int_{|z|> 1}\big[\Phi_n\big(X_{s-}+1_{[0,\sigma(X_{s-},z)]}(r)z\big)-\Phi_n(X_{s-})\big] \cN(\dif z\times\dif r\times \dif s)\\
&\quad+\int_0^t\!\!\!\int_0^{\infty}\!\!\!\!\int_{|z|\leq 1}\big[\Phi_n\big(X_{s-}+1_{[0,\sigma(X_{s-},z)]}(r)z\big)-\Phi_n(X_{s-})\big]\tilde \cN(\dif z\times\dif r\times \dif s)\\
&=\Phi_n(x)+\int_0^t\!\lambda u_n(X_s)\dif s+\int_0^t\!\Big(\big(b-b_n\big)\cdot\nabla u_n+\big(b-b_n\big)\Big)(X_s)\dif s\\
&\quad-\int_0^t\!\!\!\int_{|z|>1}\big[u_n(X_s+z)-u_n(X_s)\big]\sigma(X_s,z)\nu(\dif z)\dif s\\
&\quad+\int_0^t\!\!\!\int_0^{\infty}\!\!\!\!\int_{|z|> 1}1_{[0,\sigma(X_{s-},z)]}(r)\big[\Phi_n(X_{s-}+z)-\Phi_n(X_{s-})\big] \cN(\dif z\times\dif r\times \dif s)\\
&\quad+\int_0^t\!\!\!\int_0^{\infty}\!\!\!\!\int_{|z|\leq 1}1_{[0,\sigma(X_{s-},z)]}(r)\big[\Phi_n(X_{s-}+z)-\Phi_n(X_{s-})\big]\tilde \cN(\dif z\times\dif r\times \dif s)\\
&=:\Phi_n(x)+\cI_1+\cI_2+\cI_3+\cI_4.
\end{align*}
Now we are going to take limits for the above equality. First of all, it is easy to see that
$$
\lim_{n\rightarrow\infty}Y^n_t=\Phi(X_t)=Y_t.
$$
By the dominated convergence theorem and (\ref{unn}), we also have
\begin{align*}
\cI_2+\cI_3\rightarrow&-\!\int_0^t\!\!\!\int_{|z|>1}\big[u(X_s+z)-u(X_s)\big]\sigma(X_s,z)\nu(\dif z)\dif s\\
&+\int_0^t\!\!\!\int_0^{\infty}\!\!\!\!\int_{|z|> 1}1_{[0,\sigma(X_{s-},z)]}(r)\big[\Phi(X_{s-}+z)-\Phi(X_{s-})\big] \cN(\dif z\times\dif r\times \dif s).
\end{align*}
As for $\cI_4$, it follows from (\ref{upd}), (\ref{unn}) and the dominated convergence theorem that
\begin{align*}
&\mE\bigg|\int_0^t\!\!\!\int_0^{\infty}\!\!\!\int_{|z|\leq1}\!\!1_{[0,\sigma(X_{s-},z)]}(r)\Big[\Phi^n(X_{s-}+z)-\Phi^n(X_{s-})\\
&\qquad\qquad\qquad\qquad\qquad-\Phi(X_{s-}+z)+\Phi(X_{s-})\Big]\tilde{\cN}(\dif z\times\dif r\times \dif s)\bigg|^2\\
&=\mE\int_0^t\!\!\!\int_0^{\infty}\!\!\!\int_{|z|\leq1}\!\!1_{[0,\sigma(X_s,z)]}(r)\Big|\Phi^n_s(X_s+z)-\Phi^n_s(X_s)-\Phi_s(X_s+z)+\Phi_s(X_s)\Big|^2\nu(\dif z)\dif r\dif s\\
&=\mE\int_0^t\!\!\!\int_{|z|\leq1}\Big|\Phi^n_s(X_s+z)-\Phi^n_s(X_s)-\Phi_s(X_s+z)+\Phi_s(X_s)\Big|^2\sigma(X_s,z)\nu(\dif z)\dif s\\
&\rightarrow0,\quad\text{as}\,\, n\rightarrow\infty.
\end{align*}
Finally, Krylov's estimate (\ref{kry1}) yields that
$$
\mE\Bigg(\!\int_0^tb(X_s)-b_n(X_s)\dif s\Bigg)\leq C\|b-b_n\|_p\rightarrow0,
$$
which in turn implies by (\ref{unn}) that
\begin{align*}
\cI_1&\rightarrow \lambda\!\!\int_0^t u(X_s)\dif s,\quad\text{as}\,\, n\rightarrow\infty.
\end{align*}
Combing the above calculations, and noticing that $X_s=\Phi^{-1}(Y_s)$, we get the desired result.
\end{proof}

At the end of this section, we collect some properties of the new coefficients. For a function $f$ on $\mR^d$, set
$$
\cJ_zf(x):=f(x+z)-f(x).
$$
Suppose that:
\begin{enumerate}
\item [{\bf (H$\sigma'$)}] The global condition (\ref{s1}) holds true and (\ref{a1}) is satisfied for almost all $x,y\in\mR^d$ with $\zeta\in L^q(\mR^d)$, $q>d/\alpha$.

\item [{\bf (Hb$'$)}] For $\theta, p$ satisfying (\ref{index}),
$$
b\in  L^{\infty}(\mR^d)\cap\mW^{\theta}_p.
$$
\end{enumerate}
Then, we have:
\bl
Under {\bf (H$\sigma'$)}-{\bf (Hb$'$)}, there exist constants $C_1, C_2$ such that for a.e. $x,y\in\mR^d$,
\begin{align}
|\tilde b(x)-\tilde b(y)|\leq C_1|x-y|\cdot\Big(1+\zeta\big(\Phi^{-1}(x)\big)+\zeta\big(\Phi^{-1}(y)\big)\Big)\label{b}
\end{align}
and
\begin{align}
|\tilde g(x,z)-\tilde g(y,z)|\leq C_2|x-y|\cdot\Big(\cM|\nabla\cJ_zu|(\Phi^{-1}(x))+\cM|\nabla\cJ_zu|(\Phi^{-1}(y))\Big).\label{g}
\end{align}
Moreover, for any $p>1$ and $\gamma\in(1,2)$, it holds for all $f\in\mH^{\gamma}_p$ that
\begin{align}
\|\cJ_zf\|_{1,p}\leq C_{p,d,\gamma}|z|^{\gamma-1}\|f\|_{\gamma,p},\label{jz}
\end{align}
where $C_{p,d,\gamma}$ is a positive constant.
\el
\begin{proof}
Recall the definition of $\tilde b$ and $\tilde g$ in Lemma \ref{zvon}. Since $\sigma$ is bounded and thanks to (\ref{es2}), (\ref{upd}), (\ref{a1}), we get
\begin{align*}
|\tilde b(x)-\tilde b(y)|&\leq \lambda\big|u\big(\Phi^{-1}(x)\big)-u\big(\Phi^{-1}(y)\big)\big|+\int_{|z|>1}\!\big|u\big(\Phi^{-1}(x)+z\big)-u\big(\Phi^{-1}(y)+z\big)\big|\nu(\dif z)\\
&\quad+\int_{|z|>1}\!\big|u\big(\Phi^{-1}(x)\big)-u\big(\Phi^{-1}(y)\big)\big|\nu(\dif z)+\int_{|z|>1}\!\big|\sigma\big(\Phi^{-1}(x),z\big)-\sigma\big(\Phi^{-1}(y),z\big)\big|\nu(\dif z)\\
&\leq C_\lambda|x-y|+C_0|x-y|\Big(\zeta\big(\Phi^{-1}(x)\big)+\zeta\big(\Phi^{-1}(y)\big)\Big),
\end{align*}
which gives (\ref{b}). By (\ref{w11}), further have
\begin{align*}
|\tilde g(x,z)-\tilde g(y,z)|&=\big|\Phi\big(\Phi^{-1}(x)+z\big)-\Phi\big(\Phi^{-1}(x)\big)-\Phi\big(\Phi^{-1}(y)+z\big)+\Phi\big(\Phi^{-1}(y)\big)\big|\\
&=\big|u\big(\Phi^{-1}(x)+z\big)-u\big(\Phi^{-1}(x)\big)-u\big(\Phi^{-1}(y)+z\big)+u\big(\Phi^{-1}(y)\big)\big|\\
&=\big|\big(\cJ_zu\big)(\Phi^{-1}(x))-\big(\cJ_zu\big)(\Phi^{-1}(y))\big|\\
&\leq C_2|x-y|\cdot\Big(\cM|\nabla\cJ_zu|(\Phi^{-1}(x))+\cM|\nabla\cJ_zu|(\Phi^{-1}(y))\Big).
\end{align*}
As for (\ref{jz}), it was proved by \cite[Lemma 2.3]{Zh00}.
\end{proof}

\section{Proof of the main result}

Now, we are ready to give the proof of our main result. Comparing with the usual SDEs driven by Brownian motion or pure jump L\'evy processes, for SDEs of the form (\ref{sde2}) or (\ref{sde3}), some new tricks are needed to handle the term $1_{[0,\sigma(X_s,z)]}(r)$, as we shall see below. The point is that we have to use $L^1$-estimate as well as $L^2$-estimate to deduce the pathwise uniqueness.
\begin{proof}[Proof of Theorem \ref{main}]
The proof will consist of two steps.\\
{\bf Step 1:} We assume that {\bf (H$\sigma'$)}-{\bf (Hb$'$)} hold. It was shown in \cite[Proposition 3]{M-P} that under these conditions, there exists a martingale solution to operator $\sL$. Meanwhile, it is known that the martingale solution for $\sL$ is equivalent to the weak solution to SDE (\ref{sde2}), see \cite[Lemma 2.1]{Kurz2}. Thus, the existence and uniqueness of weak solution hold for SDE (\ref{sde2}). Thus, it suffices to show the pathwise uniqueness.

Let $X_t$ and $\hat X_t$ be two strong solutions for SDE (\ref{sde2}) both starting from $x\in\mR^d$, and set
$$
Y_t:=\Phi(X_t),\quad \hat Y_t:=\Phi(\hat X_t).
$$
Since the uniqueness if a local property, as the argument in \cite[Theorem IV. 9.1]{Ik-Wa} and \cite{Zh00}, we only need to prove by Lemma \ref{zvon} that
\begin{align}
Z_t\equiv0,\quad\forall t\geq 0, \label{66}
\end{align}
where $Z_t$ is given by
\begin{align*}
Z_{t\wedge\tau_1}&=\int_0^{t\wedge\tau_1}\!\big[\tilde b(Y_s)-\tilde b(\hat Y_s)\big]\dif s +\!\int_0^{t\wedge\tau_1}\!\!\!\!\int_0^{\infty}\!\!\!\!\int_{|z|\leq 1}\Big[\tilde g(Y_{s-},z)1_{[0,\tilde\sigma(Y_{s-},z)]}(r)\\
&\quad\qquad\quad\qquad\quad\qquad\quad-\tilde g(\hat Y_{s-}, z)1_{[0,\tilde\sigma(\hat Y_{s-},z)]}(r)\Big]\tilde \cN(\dif z\times\dif r\times \dif s)=:\cJ^{t\wedge\tau_1}_1+\cJ^{t\wedge\tau_1}_2.
\end{align*}
Set
$$
A_1(t):=\int_0^t\Big(1+\zeta(X_s)+\zeta(\hat X_s)\Big)\dif s,
$$
then following by an approximation argument as in \cite{Zh3,Zh00}, it is easy to see by (\ref{b}) that for almost all $\omega$ and every stopping time $\eta$,
\begin{align*}
\sup_{t\in[0,\eta]}\left|\cJ^{t\wedge\tau_1}_1\right|\leq C_1\!\int_0^{\tau_1\wedge\eta}|Z_s|\cdot\Big(1+\zeta(X_s)+\zeta(\hat X_s)\Big)\dif s=C_1\!\int_0^{\tau_1\wedge\eta}|Z_s| \dif A_1(s).
\end{align*}
As for the second term, write
\begin{align*}
\cJ_2^{t\wedge\tau_1}&=\int_0^{t\wedge\tau_1}\!\!\!\!\int_0^{\infty}\!\!\!\!\int_{|z|\leq 1}1_{[0,\tilde\sigma(Y_{s-},z)\wedge\tilde\sigma(\hat Y_{s-},z)]}(r)\Big[\tilde g(Y_{s-},z)1_{[0,\tilde\sigma(Y_{s-},z)]}(r)\\
&\qquad\qquad\qquad\qquad\qquad\qquad\qquad-\tilde g(\hat Y_{s-}, z)1_{[0,\tilde\sigma(\hat Y_{s-},z)]}(r)\Big]\tilde \cN(\dif z\times\dif r\times \dif s)\\
&\quad+\int_0^{t\wedge\tau_1}\!\!\!\!\int_0^{\infty}\!\!\!\!\int_{|z|\leq 1}1_{[\tilde\sigma(Y_{s-},z)\vee\tilde\sigma(\hat Y_{s-},z),\infty]}(r)\Big[\tilde g(Y_{s-},z)1_{[0,\tilde\sigma(Y_{s-},z)]}(r)\\
&\qquad\qquad\qquad\qquad\qquad\qquad\qquad-\tilde g(\hat Y_{s-}, z)1_{[0,\tilde\sigma(\hat Y_{s-},z)]}(r)\Big]\tilde \cN(\dif z\times\dif r\times \dif s)\\
&\quad+\int_0^{t\wedge\tau_1}\!\!\!\!\int_0^{\infty}\!\!\!\!\int_{|z|\leq 1}1_{[\tilde\sigma(Y_{s-},z)\wedge\tilde\sigma(\hat Y_{s-},z),\tilde\sigma(Y_{s-},z)\vee\tilde\sigma(\hat Y_{s-},z)]}(r)\Big[\tilde g(Y_{s-},z)1_{[0,\tilde\sigma(Y_{s-},z)]}(r)\\
&\qquad\qquad\qquad\qquad\qquad\qquad\qquad-\tilde g(\hat Y_{s-}, z)1_{[0,\tilde\sigma(\hat Y_{s-},z)]}(r)\Big]\tilde \cN(\dif z\times\dif r\times \dif s)\\
&=:\cJ_{21}^{t\wedge\tau_1}+\cJ_{22}^{t\wedge\tau_1}+\cJ_{23}^{t\wedge\tau_1}.
\end{align*}
We proceed to estimate each component. First, for $\cJ_{21}^{t\wedge\tau_1}$, we use the Doob's $L^2$-maximal inequality to deduce that for any stopping time $\eta$,
\begin{align*}
\mE\left[\sup_{t\in[0,\eta]}|\cJ_{21}^{t\wedge\tau_1}|\right]&\leq \mE\Bigg(\int_0^{\tau_1\wedge\eta}\!\!\!\!\int_0^{\infty}\!\!\!\!\int_{|z|\leq 1}1_{[0,\tilde\sigma(Y_s,z)\wedge\tilde\sigma(\hat Y_s,z)]}(r)\big|\tilde g(Y_s,z)-\tilde g(\hat Y_s, z)\big|^2\dif r\nu(\dif z) \dif s\Bigg)^{\frac{1}{2}}\\
&=\mE\Bigg(\int_0^{\tau_1\wedge\eta}\!\!\!\!\int_{|z|\leq 1}\big[\tilde\sigma(Y_s,z)\wedge\tilde\sigma(\hat Y_s,z)\big]\cdot\big|\tilde g(Y_s,z)-\tilde g(\hat Y_s, z)\big|^2\nu(\dif z) \dif s\Bigg)^{\frac{1}{2}}.
\end{align*}
Thus, if we set
$$
A_2(t):=\int_0^t\!\!\!\int_{|z|\leq 1}\!\Big(\cM|\nabla\cJ_zu|(X_s)+\cM|\nabla\cJ_zu|(\hat X_s)\Big)^2\nu(\dif z) \dif s,
$$
we further have by the fact that $\tilde\sigma$ is bounded and (\ref{g}) that
\begin{align*}
\mE\left[\sup_{t\in[0,\eta]}|\cJ_{21}^{t\wedge\tau_1}|\right]&\leq C_2\mE\Bigg(\int_0^{\tau_1\wedge\eta}|Z_s|^2\!\int_{|z|\leq 1}\!\Big(\cM|\nabla\cJ_zu|(X_s)+\cM|\nabla\cJ_zu|(\hat X_s)\Big)^2\nu(\dif z) \dif s\Bigg)^{\frac{1}{2}}\\
&=C_2\mE\Bigg(\int_0^{\tau_1\wedge\eta}|Z_s|^2\dif A_2(s)\Bigg)^{\frac{1}{2}}.
\end{align*}
Next, it is easy to see that for any $t\geq 0$,
$$
\cJ_{22}^{t\wedge\tau_1}\equiv0.
$$
Finally, we use the $L^1$-estimate (see \cite[P$_{174}$]{Kurz3} or \cite[P$_{157}$]{Kurz2}) to control the third term by
\begin{align*}
\mE\left[\sup_{t\in[0,\eta]}|\cJ_{23}^{t\wedge\tau_1}|\right]&\leq 2\mE\!\int_0^{\tau_1\wedge\eta}\!\!\!\!\int_0^{\infty}\!\!\!\!\int_{|z|\leq 1}1_{[\tilde\sigma(Y_s,z)\wedge\tilde\sigma(\hat Y_s,z),\tilde\sigma(Y_s,z)\vee\tilde\sigma(\hat Y_s,z)]}(r)\\
&\qquad\qquad\qquad\times\big|\tilde g(Y_s,z)1_{[0,\tilde\sigma(Y_s,z)]}(r)-\tilde g(\hat Y_s, z)1_{[0,\tilde\sigma(\hat Y_s,z)]}(r)\big|\nu(\dif z)\dif r \dif s\\
&\leq2\mE\!\int_0^{\tau_1\wedge\eta}\!\!\!\int_{|z|\leq 1}|\tilde\sigma(Y_s,z)-\tilde\sigma(\hat Y_s,z)|\cdot\Big(|\tilde g(Y_s,z)|+|\tilde g(\hat Y_s,z)|\Big)\nu(\dif z)\dif s.
\end{align*}
Since
$$
|\tilde g(x,z)|=\big|\Phi\big(\Phi^{-1}(x)+z\big)-\Phi\big(\Phi^{-1}(x)\big)\big|\leq \frac{3}{2}|z|,
$$
and taking into account of (\ref{a1}), we get
\begin{align*}
\mE\left[\sup_{t\in[0,\eta]}|\cJ_{23}^{t\wedge\tau_1}|\right]&\leq C_3\mE\!\int_0^{\tau_1\wedge\eta}\!\!\!\int_{|z|\leq 1}|\tilde\sigma(Y_s,z)-\tilde\sigma(\hat Y_s,z)|\cdot|z|\nu(\dif z)\dif s\\
&\leq C_3\mE\!\int_0^{\tau_1\wedge\eta}|Z_s|\Big(\zeta(X_s)+\zeta(\hat X_s)\Big)\dif s\leq C_3\mE\!\int_0^{\tau_1\wedge\eta}|Z_s|\dif A_1(s).
\end{align*}
Combing the above computations, and set
$$
A(t):=A_1(t)+A_2(t),
$$
we arrive at that for any stopping time $\eta$, there exists a constant $C_0$ such that
\begin{align}
\mE\left[\sup_{t\in[0,\eta]}|Z_{t\wedge\tau_1}|\right]\leq C_0\mE\!\int_0^{\tau_1\wedge\eta}|Z_s|\dif A(s)+C_0\mE\Bigg(\int_0^{\tau_1\wedge\eta}|Z_s|^2\dif A(s)\Bigg)^{\frac{1}{2}}. \label{ee}
\end{align}
By our assumption that $\zeta\in L^{q}(\mR^d)$ with $q>d/\alpha$ and the Krylov estimate (\ref{kry1}), we find that
$$
\mE A_1(t)\leq t+C\|\zeta\|_q<\infty.
$$
Meanwhile, since $p>2d/\alpha$, using the Fubini's theorem, Krylov estimate, Minkovski's inequality and taken into account of (\ref{mf}), we can get
\begin{align*}
\mE A_2(t)&=\int_{|z|\leq 1}\!\mE\!\int_0^t\Big(\cM|\nabla\cJ_zu|(X_s)+\cM|\nabla\cJ_zu|(\hat X_s)\Big)^2\dif s\nu(\dif z)\\
&\leq C\!\int_{|z|\leq 1}\!\|(\cM|\nabla\cJ_zu|)^2\|_{p/2}\nu(\dif z)= C\!\int_{|z|\leq 1}\!\|\cM|\nabla\cJ_zu|\|_{p}^2\nu(\dif z)\\
&\leq C\!\int_{|z|\leq 1}\|\nabla\cJ_zu\|_{p}^2\nu(\dif z)\leq C\!\int_{|z|\leq 1}\|\cJ_zu\|_{1,p}^2\nu(\dif z).
\end{align*}
Recall our assumption that $\theta>1-\tfrac{\alpha}{2}$. Hence, we can choose a $\gamma<\alpha$ such that
$$
2(\gamma+\theta-1)>\alpha.
$$
Consequently, it follows from (\ref{jz}) and (\ref{nu}) that
\begin{align*}
\mE A_2(t)\leq C\|u\|_{\gamma+\theta,p}^2\!\int_{|z|\leq 1}|z|^{2(\gamma+\theta-1)}\nu(\dif z)<\infty.
\end{align*}
Therefore, $t\mapsto A(t)$ is a continuous strictly increasing process. Define for $t\geq 0$ the stopping time
$$
\eta_t:=\inf\{s\geq 0: A(s)\geq t\}.
$$
Then, it is clear that $\eta_t$ is the inverse of $t\mapsto A(t)$. Since $A(t)\geq t$, we further have $\eta_t\leq t$. Taking $\eta_t$ into (\ref{ee}), we have by the change of variable
\begin{align*}
&\mE\left[\sup_{s\in[0,t]}|Z_{\tau_1\wedge\eta_s}|\right]=\mE\left[\sup_{s\in[0,\eta_t]}|Z_{s\wedge\tau_1}|\right]\\
&\leq C_0\mE\!\int_0^{\tau_1\wedge\eta_t}|Z_s|\dif A(s)+C_0\mE\Bigg(\int_0^{\tau_1\wedge\eta_t}|Z_s|^2\dif A(s)\Bigg)^{\frac{1}{2}}\\
&\leq C_0\mE\!\int_0^{\eta_t}|Z_{\tau_1\wedge s}|\dif A(s)+C_0\mE\Bigg(\int_0^{\eta_t}|Z_{\tau_1\wedge s}|^2\dif A(s)\Bigg)^{\frac{1}{2}}\\
&=C_0\mE\!\int_0^{t}|Z_{\tau_1\wedge\eta_s}|\dif s+C_0\mE\Bigg(\int_0^{t}|Z_{\tau_1\wedge\eta_s}|^2\dif s\Bigg)^{\frac{1}{2}}\leq C_0\big(t+\sqrt{t}\big)\mE\left[\sup_{s\in[0,t]}|Z_{\tau_1\wedge\eta_s}|\right].
\end{align*}
Now, take $t_0$ small enough such that
$$
C_0\big(t_0+\sqrt{t_0}\big)<1,
$$
it holds that for almost all $\omega$,
$$
\sup_{s\in[0,\eta_{t_{0}}]}|Z_{s\wedge\tau_1}|=\sup_{s\in[0,t_0]}|Z_{\tau_1\wedge\eta_s}|=0.
$$
In particular,
$$
Z_{\eta_{t_0}\wedge\tau_1}=0,\quad a.s..
$$
The same arguments as above and write for $t>t_0$,
\begin{align*}
Z_{t\wedge\tau_1}&=\int_{\eta_{t_0}\wedge\tau_1}^{t\wedge\tau_1}\!\big[\tilde b(Y_s)-\tilde b(\hat Y_s)\big]\dif s+\int_{\eta_{t_0}\wedge\tau_1}^{t\wedge\tau_1}\!\!\int_0^{\infty}\!\!\!\!\int_{|z|\leq 1}\Big[\tilde g(Y_{s-},z)1_{[0,\tilde\sigma(Y_{s-},z)]}(r)\\
&\quad\qquad\quad\qquad\quad\qquad\quad\qquad\quad\qquad-\tilde g(\hat Y_{s-}, z)1_{[0,\tilde\sigma(\hat Y_{s-},z)]}(r)\Big]\tilde \cN(\dif z\times\dif r\times \dif s),
\end{align*}
we can get
\begin{align*}
&\mE\left[\sup_{s\in[t_0,t]}|Z_{\tau_1\wedge\eta_s}|\right]=\mE\left[\sup_{s\in[\eta_{t_0},\eta_t]}|Z_{s\wedge\tau_1}|\right]\\
&\leq C_0\mE\!\int_{\eta_{t_0}\wedge\tau_1}^{\tau_1\wedge\eta_t}|Z_s|\dif A(s)+C_0\mE\Bigg(\int_{\eta_{t_0}\wedge\tau_1}^{\tau_1\wedge\eta_t}|Z_s|^2\dif A(s)\Bigg)^{\frac{1}{2}}\\
&\leq C_0\mE\!\int_{t_0}^{t}|Z_{\tau_1\wedge\eta_s}|\dif s+C_0\mE\Bigg(\int_{t_0}^{t}|Z_{\tau_1\wedge\eta_s}|^2\dif s\Bigg)^{\frac{1}{2}}\\
&\leq C_0\Big[(t-t_0)+\sqrt{(t-t_0)}\Big]\mE\left[\sup_{s\in[t_0,t]}|Z_{\tau_1\wedge\eta_s}|\right].
\end{align*}
Hence, for for almost all $\omega$,
$$
\sup_{s\in[0,\eta_{2t_{0}}]}|Z_{s\wedge\tau_1}|=0.
$$
Repeating the above arguments, we may get for any $k>0$,
$$
\sup_{s\in[0,\eta_{kt_{0}}]}|Z_{s\wedge\tau_1}|=0.
$$
Noticing that $\eta_t$ is also strictly increasing, we have for all $t\geq 0$,
$$
Z_{t\wedge\tau_1}=0,\quad a.s..
$$
Thus, (\ref{66}) is proven.\\
{\bf Step 2:} Assume now that $\sigma$ and $b$ satisfy {\bf (H$\sigma$)}-{\bf (Hb)}. For each $n\in\mN$, let $\chi_n(x)\in[0,1]$ be a nonnegative smooth function in $\mR^d$ with $\chi_n(x)=1$ for all $x\in B_n$ and $\chi_n(x)=0$ for all $x\notin B_{n+1}$. Let
$$
b_n(x):=\chi_n(x)b(x),\quad \zeta_n(x):=\chi_{n+1}(x)\zeta(x)+\cM|\nabla\chi_{n+1}|(x),
$$
and
$$
\sigma_n(x,z):=1+\chi_{n+1}(x)\sigma(x,z)\big(|z|\wedge 1\big)+\Big(1-\chi_n(x)\big(|z|\wedge 1\big)\Big)\left(1+\sup_{x\in B_{n+2}}|\sigma(x,z)|\right)\mI_{d\times d}.
$$
Then, one can check easily that (\ref{s1}) holds and $b_n$ satisfies {\bf (Hb$'$)}. Meanwhile,
\begin{align*}
&\int_{\mR^d}|\sigma_n(x,z)-\sigma_n(y,z)|(|z|\wedge1)\nu(\dif z)\leq C\!\!\int_{\mR^d}|\chi_{n+1}(x)-\chi_{n+1}(y)|\big(|z|^2\wedge1\big)\nu(\dif z)\\
&\quad\leq C\!\!\int_{\mR^d}|\sigma(x,z)-\sigma(y,z)|(|z|\wedge1)\nu(\dif z)\Big(\chi_{n+1}(x)\wedge\chi_{n+1}(y)\Big)\\
&\quad\leq C|x-y|\Big(\cM|\nabla\chi_{n+1}|(x)+\cM|\nabla\chi_{n+1}|(y)\Big)+C|x-y|\Big(\chi_{n+1}(x)\zeta(x)+\chi_{n+1}(y)\zeta(y)\Big).
\end{align*}
It is obvious that $\zeta_n\in L^q(\mR^d)$ with $q>d/\alpha$, hence {\bf (H$\sigma'$)} is also true. Therefore, for each $x\in\mR^d$, there exist a unique strong solution $X_t^n(x)$ to SDE (\ref{sde2}) with coefficients $\sigma_n$ and $b_n$. For $n\geq k$, define
$$
\varsigma_{n,k}(x):=\inf\{t\geq0: |X_t^n(x)|\geq k\}\wedge n.
$$
By the uniqueness of the strong solution, we have
$$
\mP\Big(X_t^n(x)=X_t^k(x),\,\forall t\in[0,\varsigma_{n,k}(x))\Big)=1,
$$
which implies that for $n\geq k$,
$$
\varsigma_{k,k}(x)\leq \varsigma_{n,k}(x)\leq \varsigma_{n,n}(x),\quad a.s..
$$
Hence, if we let $\varsigma_k(x):=\varsigma_{k,k}(x)$, then $\varsigma_k(x)$ is an increasing sequence of stopping times and for $n\geq k$,
$$
\mP\Big(X_t^n(x)=X_t^k(x),\,\forall t\in[0,\varsigma_{k}(x))\Big)=1.
$$
Now, for each $k\in\mN$, we can define $X_t(x):=X^k_t(x)$ for $t<\varsigma_k(x)$ and $\varsigma(x):=\lim_{k\rightarrow\infty}\varsigma_k(x)$. It is clear that $X_t(x)$ is the unique strong solution to SDE (\ref{sde2}) up to the explosion time $\varsigma(x)$ and (\ref{xt}) holds.
\end{proof}

\end{document}